\documentclass[a4paper,twoside,10pt]{article}
\usepackage{a4wide}
\usepackage{url}
\usepackage[hyperindex,pdfmark]{hyperref}

\usepackage{amsmath}	
\usepackage{bbm}
\usepackage{amssymb}
\usepackage{amsthm}
\usepackage{epic}
\usepackage{eepic}
\usepackage{subfigure}
\newcounter{pfenumi}
\newenvironment{proofenum}{\begin{list}{\it(\roman{pfenumi})}{\usecounter{pfenumi}\setlength{\leftmargin}{0pt}\setlength{\labelwidth}{-1ex}}}{\end{list}}

\def\ignore#1{}
\newtheorem{thm}{Theorem}
\newtheorem{lem}[thm]{Lemma}
\newtheorem{prop}[thm]{Proposition}
\newtheorem{cor}[thm]{Corollary}
\newtheorem{problem}{Question}
\newtheorem{que}[problem]{Question}
\newtheorem*{claim}{Claim}

\theoremstyle{definition}
\newtheorem*{defn}{Definition}
\newtheorem*{notation}{Notation}
\newtheorem*{rmk}{Remark}

\allowdisplaybreaks[2]
\DeclareMathOperator{\adj}{\operatorname{adj}}
\DeclareMathOperator{\nbd}{\operatorname{nbd_\delta}}
\DeclareMathOperator{\tr}{\operatorname{tr}}
\DeclareMathOperator{\kmf}{{\mathcal K}_\delta}
\DeclareMathOperator{\nmf}{{\mathcal N}_\delta}
\DeclareMathOperator{\lkmf}{\mathcal{LK}_\delta}
\DeclareMathOperator{\lnmf}{\mathcal{LN}_\delta}
\DeclareMathOperator{\grad}{\nabla\!}
\DeclareMathOperator{\supp}{{\operatorname{supp}}}
\DeclareMathOperator{\dist}{{\operatorname{dist}}}

\def\trns{\mathrm{t}} 
\def\ball{{\mathbb B}^{n-1}}

\def\eps{\varepsilon}
\def\A{{\mathcal A}}
\def\Abb{{\mathbb A}}
\def\B{{\mathcal B}}
\def\Bbb{{\mathbb B}}
\def\Complex{{\mathbb C}}
\def\d{\,{\mathrm d}}
\def\G{{\mathcal G}}
\def\H{{\mathbb H}}
\def\Lcal{{\mathcal L}}
\def\Ncal{{\mathcal N}}
\def\Q{{\mathbb Q}}
\def\R{{\mathbb R}}
\def\Z{{\mathbb Z}}

\def\1{{\mathbbm 1}}

\begin{document}
\title{Kakeya sets of curves}
\author{Laura Wisewell%
\thanks{This work formed part of my PhD thesis \protect\cite{wisewell:thesis}.  The encouragement and help of my supervisor Professor Tony Carbery, and the financial support of the EPSRC and the Seggie-Brown Trust are gratefully acknowledged.\newline2000 Mathematics Subject Classification 42B25}}
\maketitle
\section{Introduction}
\begin{figure}[!hbt]
\centering
\setlength{\unitlength}{0.00083333in}
\begin{picture}(4148,2181)(0,-10)
\path(2055,324)(2175,354)(2055,384)
\path(2175,354)(1875,354)
\path(2775,1254)(3075,1254)(3075,954)
	 (2775,954)(2775,1254)
\texture{8101010 10000000 444444 44000000 11101 11000000 444444 44000000 
	 101010 10000000 444444 44000000 10101 1000000 444444 44000000 
	 101010 10000000 444444 44000000 11101 11000000 444444 44000000 
	 101010 10000000 444444 44000000 10101 1000000 444444 44000000 }
\shade\path(2775,1854)(3075,1854)(3075,1554)
	 (2775,1554)(2775,1854)
\path(375,354)(375,2154)
\path(405,2034)(375,2154)(345,2034)
\path(975,954)(2175,354)(2175,354)
	 (375,354)(375,954)(975,954)
\shade\path(375,354)(375,804)(1275,804)
	 (2175,354)(375,354)
\dashline{60}(375,354)(2175,2154)
\put(1200,60){$\frac 12$}
\put(2175,60){$\frac 1r$}
\put(75,2000){$\frac1s$}

\put(3225,1550){known}
\put(3225,950){conjectured}
\put(1575,700){$s=\frac{(n+1)r'}{n-1}$}
\put(-100,650){$\frac{n-1}{2(n+1)}$}
\put(-100,950){$\frac{n-1}{2n}$}
\end{picture}
\caption{Exponents for H\"ormander's conjecture}\label{fig:horpq}
\end{figure}
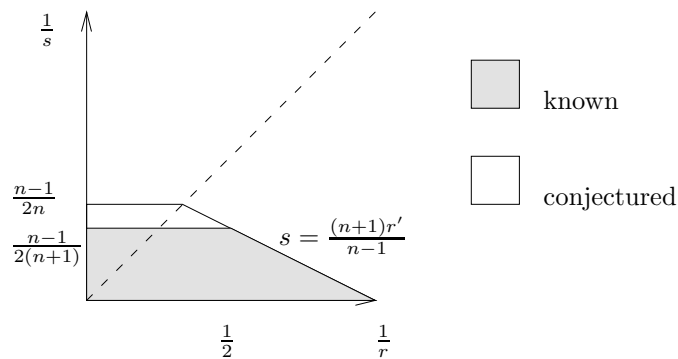

In this paper we investigate an analogue for curves of the famous
Kakeya conjecture about straight lines.  The simplest version of the
latter asks whether a set in $\R^n$ that includes a unit line segment
in every direction must necessarily have dimension $n$.  The analogue
we have in mind replaces the line segments by curved arcs from a
specified family.  (This is a quite different problem from that
considered by Minicozzi and Sogge \cite{minicozzi:negative} who looked
at geodesics in curved space.)  The families of curves we are
interested in arise from H\"ormander's conjecture in harmonic
analysis, which deals with oscillatory integral operators of the form
\begin{equation}\label{TN}
T_Nf(x):=\int_{\R^{n-1}}e^{iN\varphi (x,y)}a(x,y)f(y)\d
y.\end{equation}
Here $x\in\R^n$, $y\in\R^{n-1}$, $a$ is some smooth cut-off, and
the phase function $\varphi$ is assumed to be smooth on the support of
$a$ and to have the following properties:
\begin{gather}
\mbox{The matrix }\frac{\partial^2\varphi}{\partial x\partial y}(x,y)\mbox{ has full rank }n-1.\label{rank}\\
{\mbox{For all }\theta\in S^{n-1}\mbox{ the map }}y\mapsto\theta\cdot\frac{\partial\varphi}{\partial x}(x,y){\mbox{ has only non-degenerate critical points.}}\label{nondeg}
\end{gather}
In \cite{bourgain:lp} it is observed that by making appropriate changes of variable, any phase satisfying these criteria can be expressed in the form
\begin{equation}\label{generalphase}
\varphi(x,y)=y^\trns x'+x_n y^\trns Ay+O(|x_n||y|^3+|x|^2|y|^2)
\end{equation}
with $A$ an invertible $(n-1)\times(n-1)$ matrix and $\trns$ denoting transpose.  H\"ormander showed
that both Restriction and Bochner-Riesz problems can be formulated as
special cases of operators $T_N$, which prompted him to ask the
following:
\begin{que}[H\"ormander \cite{hormander:oscill}]\label{horconj}  Is it true that for every $\varphi$ satisfying the above properties, the operator $T_N$ has the bound 
\begin{equation}
\label{Tineq}
\|T_N f\|_s\lesssim N^{-n/s}\|f\|_r\ ?\end{equation}
for $\frac 1s<\frac{n-1}{2n}$ and $\frac 1s\leq\frac{n-1}{n+1}\frac {1}{r'}$?
\end{que}
H\"ormander himself proved this for $n=2$ \cite{hormander:oscill}, and
in higher dimensions it has been proved for $s\geq \frac{2(n+1)}{n-1}$
by Stein \cite{stein:oscill}.  The known and conjectured regions are
shown in Figure~\ref{fig:horpq}.

 It was a great surprise in 1991 when Bourgain \cite{bourgain:lp}
 disproved H\"or\-man\-der's conjecture.  Roughly, he showed that in
 dimension three for most phases the best exponent $s$ is strictly
 greater than $\frac{2n}{n-1}=3$, and that there exist phases where
 the known value $\frac{2(n+1)}{n-1}=4$ is the best.  More precisely:
\begin{thm}[Worst Case, \cite{bourgain:lp}] 
\label{thm:worst}
 In dimension three there is a phase function, namely
 $$\varphi(x,y)=x_1y_1+x_2y_2+x_3y_1y_2+\frac12x_3^2y_1^2$$
 for which \eqref{Tineq} fails for all $s<4$, even with $r=\infty$.
\end{thm}
\begin{thm}[Generic Failure, \cite{bourgain:lp}]
\label{thm:genfail}
 In dimension three, if $\varphi$ has the property that
\begin{equation}
\label{genericfailure}
\left.\frac{\partial^2}{\partial y^2}\left(\frac{\partial^2\varphi}{\partial x_3^2}\right)\right|_{x=0,y=0}\textrm{is not a multiple of }\left.\frac{\partial^2}{\partial y^2}\left(\frac{\partial\varphi}{\partial x_3}\right)\right|_{x=0,y=0}
\end{equation}
 then the inequality \eqref{Tineq} cannot hold even for $r=\infty$ unless $s\geq 118/39 >3$.
\end{thm}
His method was to link the oscillatory integral problems to Kakeya-type problems about curves.  Diagrammatically we have the following chain of implications:
\begin{equation}
\label{implications}
\begin{array}{c}\mbox{Oscillatory}\\\mbox{integral estimates}\end{array}\implies
\begin{array}{c}\mbox{Kakeya maximal}\\\mbox{function estimates}\end{array}\implies
\begin{array}{c}\mbox{Kakeya sets}\\\mbox{have large dimension.}\end{array}
\end{equation}
(Compare this with Fefferman's famous counterexample to the disc conjecture \cite{feffC:ball}, which used the fact that Kakeya sets can have measure zero to disprove a Bochner-Riesz type estimate.)

Both Restriction and Bochner-Riesz give rise to straight line problems ---strictly speaking, Restriction implies the Kakeya conjecture while Bochner-Riesz implies the Nikodym conjecture, which is like the Kakeya conjecture but with the roles of positions and directions exchanged.  In the straight line case these are equivalent and so attention has been focused entirely on the former, but for any fixed class of curves they are different, as will be seen in Section~\ref{sec:neg} where we consider quadratic curves.  The relationship between the two is explored in \cite{carbery:paraboloid}, \cite{tao:brimplies} and \cite{carberyw:3notes}.

Since 1991 there has been much progress on the straight line Kakeya problem, with contributions from Bourgain, Wolff, Katz, \L aba, Tao and Schlag.  Our aim in this paper is to apply some of these new techniques to the curved case.  As we shall see, we can prove positive results for certain families of curves, which may indicate that their corresponding phase functions allow reasonably good non-trivial bounds for the operators $T_N$.

We begin by giving precise definitions and brief proofs of the implications above.  Then in Section~\ref{sec:trivial} we prove the so-called ``trivial bound'' (a maximal function estimate implying that the sets have dimension at least $\frac{n+1}{2}$), which holds for a very broad class of curves.  From Section~\ref{sec:neg} onwards we restrict to quadratic curves, and demonstrate (in Theorem~\ref{thm:easygenfail}) that these are still general enough to exhibit the pathological behaviour discovered by Bourgain.  We then tackle the maximal function problem by means of geometry, proving a result (Theorem~\ref{thm:katz}) corresponding to the lower bound $\frac{n+2}{2}$ for the dimension of Nikodym sets of parabolas satisfying a certain algebraic condition.  Finally we look at arithmetic methods, and obtain (in Theorems~\ref{thm:bestnik}, \ref{thm:nastynik} and \ref{thm:kak4}) lower bounds of the form $\alpha n+\beta$ with $\alpha>1/2$ for the dimension of various sets of curves, including a bound for the Nikodym sets of the previous section which equals the best currently known for straight lines.
\subsection{The relevance of Kakeya with curves}
Given a phase function $\varphi$ and cutoff $a\in C_c^\infty$ as in \eqref{TN}, define curves and
 curved tubes as follows:
\begin{notation}
 Let $y,\omega\in\ball$  and let $\delta>0$ be a thickness.  Define
\begin{align*}
\Gamma_y(\omega)&:=\left\{x\in\R^n:\grad_y\varphi(x,y)=\omega,\,(x,y)\in\supp(a)\right\}\\
 T^\delta_y(\omega)&:=\left\{x\in\R^n:|\!\grad\!_y\varphi(x,y)-\omega|<\delta,\,(x,y)\in\supp(a)\right\}
\end{align*}
 to be the curve ``centre'' $\omega$ in ``direction'' $y$ and the corresponding $\delta$-tube.
\end{notation}
Here $\ball$ denotes a ball in $\R^{n-1}$ of some constant radius: for quadratic curves the unit ball will do, but more generally we will need to choose the constant to depend on $\varphi$, although of course larger sets of directions can then be handled simply by taking unions of small enough balls. 

Using the rank condition \eqref{rank} and the implicit function theorem we see that $\Gamma_y(\omega)$ is indeed a smooth curve.
The descriptions ``centre'' and ``direction'' are to aid intuition; in some cases the meaning of the variables may in fact be the other way round.  On one hand, the Restriction problem for the paraboloid corresponds to the phase $\varphi(x,y)=y^\trns x'+x_ny^\trns y$, so that $\Gamma_y(\omega)$ is a straight line centred at $(\begin{smallmatrix}\omega\\0\end{smallmatrix})$ in direction $(\begin{smallmatrix}y\\1\end{smallmatrix})$.  But on the other hand, Bochner-Riesz for the paraboloid has the \label{brphase}phase $\varphi(x,y)=\frac{1}{x_n}y^\trns x'+\frac{1}{x_n}y^\trns y$, and $\Gamma_y(\omega)$ is still a straight line, but with centre $(\begin{smallmatrix}y\\0\end{smallmatrix})$ and direction $(\begin{smallmatrix}\omega\\1\end{smallmatrix})$.  

 Because of the smoothness of $\varphi$, the tubes have the following ``doubling property'': Suppose that $|y-\bar{y}|<\delta$ and $|\omega-\bar{\omega}|<\delta$.  Then $T^{C\delta}_y(\omega)\supset T^{\delta}_{\bar{y}}(\bar{\omega})$ for some constant $C$ depending only on $\varphi$.  This will often allow us to consider only finite $\delta$-separated collections of tubes.

 By analogy with the straight line case, we define the following sets:
\begin{defn}[Curved Kakeya set]A set $E\subset\R^n$ is a curved Kakeya set (associated to $\varphi$) if for all $y\in\ball$ there exists an $\omega\in\ball$ such that $\Gamma_y(\omega)\subset E$.  
\end{defn}
\noindent So for the Restriction phase above, this is the usual definition of a Kakeya set, while for the Bochner-Riesz phase this is a set that includes a line segment in some direction through every one of a large set of points, which might be termed a {\em Nikodym set} although this is not quite the same definition as is usually given in, say, \cite{falconer:frac}.

As we shall see shortly, curved Kakeya sets need not have full dimension, so rather than a conjecture we have a question:
\begin{que}Given a phase function, what is the minimum possible dimension for its corresponding Kakeya sets?  For which curves must the dimension be exactly $n$?
\end{que}
We can ask this about either the usual Hausdorff notion of dimension, or more weakly about the {\em Minkowski dimension}.  This is always greater than or equal to the Hausdorff and is simpler to use: a set $E$ has (upper) Minkowski dimension at least $d$ if and only if its $\delta$-neighbourhood $\nbd(E)$ has Lebesgue measure satisfying
$|\nbd(E)|\gtrsim\delta^{n-d}$.
This is the notion we shall use in Section~\ref{sec:arith} when applying the arithmetic techniques.

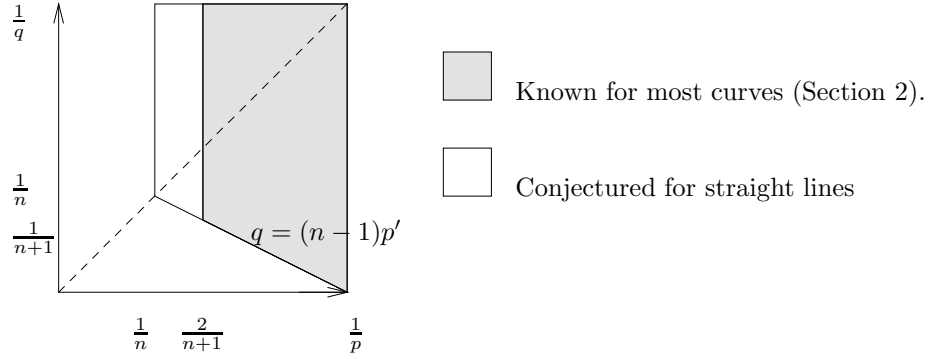
\begin{figure}[!hbt]
\centering
\setlength{\unitlength}{0.00083333in}
\bigskip
 {\renewcommand{\dashlinestretch}{30}
\begin{picture}(4103,2181)(500,-10)
\path(330,2034)(300,2154)(270,2034)
\path(300,2154)(300,354)(2100,354)
	 (900,954)(900,2154)(1200,2154)
\path(1980,324)(2100,354)(1980,384)
\path(2100,354)(1800,354)
\texture{8101010 10000000 444444 44000000 11101 11000000 444444 44000000 
	 101010 10000000 444444 44000000 10101 1000000 444444 44000000 
	 101010 10000000 444444 44000000 11101 11000000 444444 44000000 
	 101010 10000000 444444 44000000 10101 1000000 444444 44000000 }
\shade\path(1200,2154)(1200,804)(2100,354)
	 (2100,2154)(1200,2154)
\path(1200,2154)(1200,804)(2100,354)
	 (2100,2154)(1200,2154)
\dashline{60}(300,354)(2100,2154)
\path(2700,1254)(3000,1254)(3000,954)
	 (2700,954)(2700,1254)
\shade\path(2700,1854)(3000,1854)(3000,1554)
	 (2700,1554)(2700,1854)
\put(750,54){$\frac 1n$}
\put(1050,54){$\frac2{n+1}$}
\put(2100,54){$\frac1p$}
\put(0,2004){$\frac1q$}
\put(0,954){$\frac1n$}
\put(0,654){$\frac1{n+1}$}
\put(1500,680){$q=(n-1)p'$}
\put(3150,1554){Known for most curves (Section~\ref{sec:trivial}). }
\put(3150,954){Conjectured for straight lines}
\end{picture}
 }
\caption{Region where the $L^p\to L^q$ norms of the Kakeya maximal operator should be at most $\delta^{1-n/p-\eps}$}\label{fig:kakeyapq}
\end{figure}

More difficult questions about overlap of tubes can be posed in terms of
maximal operators.
\begin{defn}[Curved maximal operator]The curved Kakeya maximal function (associated to $\varphi$ and of eccentricity $1/\delta$) is the operator that takes a function $f$ on $\R^n$ to the function $\kmf f$ on $\ball$ given by
 $$\kmf f(y):=\sup_{\omega\in\ball}
\frac1{|T_y(\omega)|}\int_{T_y(\omega)}|f(x)|\d x.$$ 
\end{defn}
\noindent In the straight line case it is conjectured that this should have $L^n\to L^n$ operator norm at most $\delta^{-\eps}$.  This and the bounds which follow by interpolation are shown in Figure~\ref{fig:kakeyapq}.
However as we shall see shortly, for many families of curves the $L^n\to L^n$ bound is false, so we have another question.
\begin{que}
\label{que:max} Given a phase function, how can we find  the best bound for the corresponding maximal functions?  For which phases is an $L^n\to L^n$ bound of order $\delta^{-\eps}$ possible?
\end{que}
These estimates imply lower bounds for the dimension of the sets in the following way.
\begin{prop}[Maximal implies dimension]\label{prop:maxtohaus}
Assume that for some $\varphi$ an estimate $\|\kmf f\|_{q,\infty}\leq C\delta^{-\alpha}\|f\|_{p,1}$ holds.  Then the corresponding  Kakeya sets  have Hausdorff dimension at least $n-p\alpha$.
\end{prop}
So a sharp $L^p\to L^q$ bound  implies that the sets have dimension at least $p$.  The implication is easy if we use Minkowski dimension: simply let $f$ be the characteristic function of the $\delta$-neighbourhood of the set, and the required estimate follows.  The proof for the Hausdorff dimension is similar but requires a dyadic pigeonholing argument, following exactly that given by Wolff for the straight line case in \cite[Lemma 1.6]{wolff:recentkakeya}, with an implicit function argument to obtain suitable parameterisations of the curves. 

So we have the second of the implications \eqref{implications}.  We now turn to the first which relates the above to the oscillatory integrals $T_N$.  The proof will show why the curves we have defined are natural.
\begin{prop}\label{prop:oscillimplies}
 For some phase $\varphi$ satisfying H\"ormander's criteria \eqref{rank} and \eqref{nondeg}, suppose that $\| T_Nf\|_s\lesssim N^{-n/s}\|f\|_r$. Then the corresponding curved Kakeya maximal
 function is of restricted weak type (p,q) with norm at most $\delta^{-2(n/p-1)}$, where 
 $p=(s/2)'$ and $q=(r/2)'$,
\end{prop}

 To prove this, and also to prove the estimates for the maximal
 functions in later sections, it is helpful to linearise the maximal
 function so that instead of an $L^p$ bound we can prove a ``covering
 lemma'' similar to those in \cite{carbery:covering}.
\begin{defn}[Linearised operator] \label{def:LK}Decompose $\ball$ into disjoint $\delta$-cubes $Q_j$ for $j\in\ball\cap\delta\Z^{n-1}$.
To each index $j$ associate a tube $T_j=T^\delta_{y_j}(\omega_j)$ where $y_j\in Q_j$ and $\omega_j\in\ball$.  Define a linearisation of $\kmf$ by
 $$\lkmf f(y):=\sum_j\1_{Q_j}(y)\frac{1}{|T_j|}\int_{T_j}f(x)\d x.$$
\end{defn}
Now taking the adjoint of this operator puts the problem in the following useful form:
\begin{lem}[Covering lemma]\label{lem:covering}
 Let $\{T_j\}_{j=1}^M$ be $1\times\delta$-tubes with centres
 $\omega_j$ and directions $y_j$ (where both of these are in $\ball$).
 Then the estimate
 $$\left\|\sum_{j=1}^M\1_{T_j}\right\|_{p'}\leq
A(\delta)(\delta^{n-1}M)^{1/q'}$$ holds for all choices of $y_j\in Q_j$ with
arbitrary $\omega_j$ if and only if the (curved) Kakeya maximal function is of
weak type $(p,q)$ with constant $A(\delta)$.  
\end{lem}
This is easy to prove, and the details are given in \cite{wisewell:thesis}.

 Now we use the linearisation to sketch a proof of Proposition~\ref{prop:oscillimplies}, which will show the reason for the definition of
 the curves.  This proof is similar to that given by Wolff in
 \cite[pp.~153--154]{wolff:recentkakeya} for the Restriction problem,
 but incorporating ideas found in Bourgain's ``generic failure'' proof
 for curves \cite[pp.~326--327]{bourgain:lp}.
\begin{proof}[Proof of Proposition~\ref{prop:oscillimplies}]
 Suppose that we are given tubes $T_j$ with directions $y_j\in Q_j$ and arbitrary centres $\omega_j$ as above.  Set 
 $$f(y)=\sum_{j=1}^M\eps_j e^{-iN\omega_j.(y-y_j)}\1_{Q_j}(y)$$ where the $\eps_j$ are random signs.  Then $T_Nf$ is a sum of integrals $\int_{Q_j}e^{iN\left(\varphi(x,y)-\omega_j.(y-y_j)\right)}a(x,y)\d y$, and it is easy to see that if we choose $N\sim\delta^{-2}$ then the phase is roughly constant for $x\in T_j$ so that the integral is at most $\delta^{n-1}\1_{T_j}(x)$.  Applying Khinchin's inequality and the assumed bound for $\|T_N\|_{r\to s}$ gives the covering lemma required.
\end{proof}
An
immediate corollary is that the optimal $s=\frac{2n}{n-1}$ would
imply the best estimate $\|\kmf\|_{n\to n}\lesssim\delta^{-\eps}$ for the curved Kakeya maximal function and full dimension for the sets.  However, away from the optimal exponent this correspondence becomes very poor: Stein's result $s=\frac{2(n+1)}{n-1}$ merely implies that the sets have dimension at least $1$.
\section{``Trivial'' results for most curves}
\label{sec:trivial}
Before we go on to discuss the dimension of curved Kakeya and Nikodym sets we should first mention that they can indeed have measure zero.  This is proved in \cite{wisewell:meas0} by adapting a result due to Sawyer, and in fact applies to more general problems of surfaces lying in sets of measure zero.

Also we point out that, as with straight lines, the problem is entirely understood in dimension $2$ since H\"ormander's conjecture is true in the plane, and so the implications \eqref{implications} give the best possible bounds for the maximal functions and set dimensions.

The next simplest result is the $\frac{n+1}{2}$ bound.  For the maximal function with straight lines this was first proved in 1986 by Christ, Duoandikoetxea and
 Rubio de Francia \cite{cdr:max} using Fourier transform methods.
 Since then, more geometric proofs have been given.  
One of these is a two-slice version of the arithmetic methods which we look at in Section~\ref{sec:arith} for the set dimension problem.
Here we use the ``bush argument''
 of Bourgain \cite{bourgain:besicovitchtype} to obtain the stronger maximal function estimate.  In a sense this is analogous to the two-slice method by point-line duality, since the main idea
 is an estimate for the size of the intersection of two different
 tubes (compare with Lemma~\ref{lem:derivemat}). 
So we begin by looking at the way curved tubes can intersect, and
 prove the $L^{\frac{n+1}{2}}$ bound for a very broad class of families of curves.
We then show geometrically why the non-degeneracy
 criterion \eqref{nondeg} is crucial, by providing a counterexample to
 the maximal function result if it is not assumed.  We also mention the
 slightly curious fact that even without non-degeneracy, the set dimension
 result still holds.

At this point we restrict our attention to phases of a slightly simpler form than \eqref{generalphase}, namely those in which the higher-order terms depend only on $x_n$ and not on $x'$.  For convenience in later sections we write these as
\begin{equation}\label{xnonlyphase}
\varphi(x,y)=y^\trns M(x_n)x'+\tilde\varphi(x_n,y),
\end{equation}
where $M:\R\to GL(n-1,\R)$ is a matrix-valued
 function.  The curves corresponding to this can be parametrised by the height $x_n$ as follows:
\begin{equation}\label{xncurve}
\Gamma_y(\omega)=\left\{\begin{pmatrix}M(x_n)^{-1}\left[\omega-\grad_y\tilde\varphi(x_n,y)\right]\\x_n\end{pmatrix}:(x,y)\in\supp(a)\right\},
\end{equation}
where the matrix inverse exists by the rank condition \eqref{rank}. This notation is introduced because the phases we want to look at in Sections~\ref{sec:neg}--\ref{sec:arith} which give rise to parabolic curves are more conveniently expressed in the form \eqref{xnonlyphase} than as in \eqref{generalphase} where we had $M(x_n)\equiv I$.

We also now specify the radius of the ball $\ball$.  Given a phase $\varphi$ of the form \eqref{xnonlyphase}, write $\psi=\frac{\partial}{\partial x_n}\grad_y\tilde\varphi$.  H\"ormander's criterion \eqref{nondeg} tells us that the matrix $\frac{\partial}{\partial y}\psi$ has non-zero determinant throughout the support of the cutoff $a$, so let $k$ be the minimum absolute value of its eigenvalues on this support.  Then by the definition of the derivative, find $\rho>0$ such that
\begin{equation}\label{rho}
\frac{|\psi (x_n,y)-\psi (z_n,\bar{y})-D\psi (x_n,\bar{y})(y-\bar{y})|}{|y-\bar{y}|}<k/2\end{equation}
whenever $|y-\bar{y}|\leq \rho$.  This constant $\rho$ depends only on $\varphi$, and $\ball$ will be taken to mean the ball of this radius from now on.  (In the case of quadratic phases which we consider in the next sections, the above fraction is zero and so this issue does not arise.)

We now state and prove the crucial estimate for the size of the intersection of two tubes.  Call two tubes $T_y(\omega)$, $T_{\bar y}(\bar\omega)$ {\em $d$-separated} if $|y-\bar y|\geq d$.
 \begin{lem}\label{lem:intersec}
 Assuming \eqref{nondeg}, there is a constant $C$ depending only on $\varphi$ such that if two $d$-separated $\delta$-tubes corresponding to curves of the form \eqref{xncurve} meet, then the diameter of their intersection is at most $C\delta/d$.
 \end{lem}
 \begin{proof}
Curves of the form \eqref{xncurve} have tangents given by
$$\begin{pmatrix}-M(x_n)^{-1}M'(x_n)M(x_n)^{-1}\left[\omega-\grad_y\tilde\varphi(x_n,y)\right]-M(x_n)^{-1}\frac\partial{\partial x_n}\grad_y\tilde{\varphi}(x_n,y)\\1\end{pmatrix}.$$
If two different curves $\Gamma_y(\omega)$ and $\Gamma_{\bar y}(\bar\omega)$ meet at height $x_n=t_0$, then
$$M(t_0)^{-1}\left[\omega-\grad_y\tilde\varphi(t_0,y)\right]=M(t_0)^{-1}\left[\bar\omega-\grad_{y}\tilde\varphi(t_0,\bar y)\right]
$$
 and so the  difference between their tangents is simply
\begin{align*}
\left|M(t_0)^{-1}\left[\psi (t_0,y)-\psi (t_0,\bar{y})\right]\right| &\gtrsim |D\psi (t_0,\bar{y})(y-\bar{y})|-|\psi (t_0,y)-\psi (t_0,\bar{y})-D\psi (t_0,\bar{y})(y-\bar{y})|\\
&>\frac k2|y-\bar{y}|\qquad\mbox{by \eqref{rho}.}
\end{align*}
So the tangents are at an angle comparable to $|y-\bar{y}|$ and so the diameter of the intersection is at most $\frac{\delta}{\delta+|y-\bar{y}|}$, giving the result claimed.
\end{proof}
This allows us to prove the $L^{\frac{n+1}{2}}$ bound using the bush argument.
 \begin{thm}
 \label{thm:trivial}
 Assuming \eqref{nondeg}, the curved Kakeya maximal function $\kmf$ corresponding to curves of the form \eqref{xncurve} satisfies
 \begin{equation}
 \label{estimate}
 \|\kmf f\|_q\leq C_\eps \delta^{-(n/p -1+\eps)}\|f\|_p
 \end{equation}
  for $1\leq p\leq\frac{n+1}{2}$ and $1\leq q\leq (n-1)p'$.
 \end{thm}
 \begin{proof}It is enough to prove a restricted weak type estimate at the endpoint, since this implies strong type at the cost of an additional log \cite[p.~48]{carbery:igari}.  The proof follows exactly the bush argument for the straight line case, a suitable version of which is given in \cite{wolff:recentkakeya} or \cite{wisewell:thesis}.
 \end{proof}
 So the so-called ``trivial bound'' holds for {\em all}
 curves of the form \eqref{xncurve}, and in particular, it is true for the ``worst case'' example
 of Bourgain. \label{skunk} That example had curves given by
 $$\begin{pmatrix}
\omega_1-x_ny_2-x_n^2y_1\\\omega_2-x_ny_1\\x_n\end{pmatrix}.$$ \label{pushtosurf}If we choose
$\omega_1=0,\ \omega_2=-y_2$ then we see that each curve lies in the
surface $x_1=x_2x_3$.  So the Kakeya set has dimension two, and for this $\varphi$ the ``trivial'' bound is in fact
best possible.

This suggests that $\frac{n+1}{2}$ for the set dimension and maximal
function ought to correspond to the exponent $s=\frac{2(n+1)}{n-1}$ in
H\"ormander's conjecture, since this is the result that is known to
be true for all phases and cannot be improved for Bourgain's
example.  However, the implication proved in Proposition~\ref{prop:oscillimplies} is weaker; one feels that the factor of $2$ in the
power of $\delta$ we obtained should not be there.

\begin{figure}[!hbt]
 \centering
 \subfigure[Proper Intersection]{
 \label{subfig:proper}
 \begin{minipage}[b]{0.5\textwidth}
 \centering
 \setlength{\unitlength}{0.0005in}
 {\renewcommand{\dashlinestretch}{30}
 \begin{picture}(1674,2439)(0,-10)
 \path(612,2412)(612,12)
 \path(912,12)(912,2412)
 \path(12,12)(13,14)(14,20)
	 (17,30)(22,45)(28,66)
	 (36,94)(46,127)(57,166)
	 (70,209)(83,255)(98,304)
	 (113,354)(128,405)(143,455)
	 (157,503)(172,550)(185,594)
	 (198,637)(211,677)(223,714)
	 (234,750)(245,783)(255,814)
	 (265,844)(275,872)(284,899)
	 (294,925)(303,950)(312,974)
	 (323,1004)(335,1033)(347,1062)
	 (359,1090)(371,1119)(383,1147)
	 (396,1175)(409,1202)(422,1229)
	 (435,1256)(449,1283)(462,1309)
	 (475,1334)(489,1359)(502,1383)
	 (515,1406)(528,1428)(541,1450)
	 (553,1470)(565,1490)(577,1509)
	 (589,1527)(601,1545)(612,1562)
	 (625,1581)(637,1599)(650,1618)
	 (663,1637)(677,1655)(691,1675)
	 (705,1694)(720,1713)(735,1733)
	 (750,1753)(765,1773)(781,1792)
	 (797,1812)(813,1831)(829,1850)
	 (844,1869)(860,1887)(875,1905)
	 (891,1923)(906,1940)(921,1957)
	 (937,1974)(951,1990)(966,2006)
	 (982,2023)(998,2040)(1015,2058)
	 (1033,2077)(1052,2097)(1073,2118)
	 (1095,2141)(1119,2165)(1144,2191)
	 (1170,2218)(1197,2246)(1225,2273)
	 (1252,2301)(1277,2326)(1300,2350)
	 (1320,2370)(1336,2386)(1348,2398)
	 (1356,2406)(1360,2410)(1362,2412)
 \path(312,12)(313,14)(314,20)
	 (317,30)(322,45)(328,66)
	 (336,94)(346,127)(357,166)
	 (370,209)(383,255)(398,304)
	 (413,354)(428,405)(443,455)
	 (457,503)(472,550)(485,594)
	 (498,637)(511,677)(523,714)
	 (534,750)(545,783)(555,814)
	 (565,844)(575,872)(584,899)
	 (594,925)(603,950)(612,974)
	 (623,1004)(635,1033)(647,1062)
	 (659,1090)(671,1119)(683,1147)
	 (696,1175)(709,1202)(722,1229)
	 (735,1256)(749,1283)(762,1309)
	 (775,1334)(789,1359)(802,1383)
	 (815,1406)(828,1428)(841,1450)
	 (853,1470)(865,1490)(877,1509)
	 (889,1527)(901,1545)(912,1562)
	 (925,1581)(937,1599)(950,1618)
	 (963,1637)(977,1655)(991,1675)
	 (1005,1694)(1020,1713)(1035,1733)
	 (1050,1753)(1065,1773)(1081,1792)
	 (1097,1812)(1113,1831)(1129,1850)
	 (1144,1869)(1160,1887)(1175,1905)
	 (1191,1923)(1206,1940)(1221,1957)
	 (1237,1974)(1251,1990)(1266,2006)
	 (1282,2023)(1298,2040)(1315,2058)
	 (1333,2077)(1352,2097)(1373,2118)
	 (1395,2141)(1419,2165)(1444,2191)
	 (1470,2218)(1497,2246)(1525,2273)
	 (1552,2301)(1577,2326)(1600,2350)
	 (1620,2370)(1636,2386)(1648,2398)
	 (1656,2406)(1660,2410)(1662,2412)
 \end{picture}
 }\end{minipage}}%
 \subfigure[Tangential Intersection]{
 \label{subfig:improper}
 \begin{minipage}[b]{0.5\textwidth}
 \centering
 \setlength{\unitlength}{0.0005in}
 {\renewcommand{\dashlinestretch}{30}
 \begin{picture}(1149,2439)(0,-10)
 \path(12,2412)(12,12)
 \path(312,12)(312,1287)(312,2412)
 \path(1137,12)(1135,14)(1130,19)
	 (1122,27)(1109,40)(1092,58)
	 (1070,80)(1044,106)(1015,135)
	 (984,167)(952,200)(919,234)
	 (886,267)(855,299)(825,330)
	 (796,360)(770,387)(745,413)
	 (722,438)(701,460)(682,482)
	 (663,502)(646,521)(630,539)
	 (614,557)(600,574)(583,595)
	 (566,615)(550,634)(535,654)
	 (520,675)(505,695)(491,715)
	 (477,736)(463,757)(451,778)
	 (438,800)(427,821)(416,842)
	 (406,863)(396,884)(388,904)
	 (380,924)(372,945)(366,965)
	 (360,984)(354,1004)(350,1024)
	 (345,1045)(341,1066)(337,1088)
	 (334,1111)(331,1134)(329,1158)
	 (327,1183)(325,1208)(324,1234)
	 (324,1261)(324,1287)(325,1313)
	 (326,1340)(328,1366)(331,1391)
	 (334,1416)(337,1440)(341,1463)
	 (346,1486)(351,1508)(356,1529)
	 (362,1549)(369,1570)(376,1590)
	 (384,1609)(392,1629)(401,1649)
	 (411,1669)(422,1690)(434,1710)
	 (446,1730)(459,1751)(472,1771)
	 (486,1791)(501,1811)(515,1831)
	 (530,1850)(545,1869)(560,1887)
	 (576,1905)(591,1923)(606,1940)
	 (621,1957)(637,1974)(651,1990)
	 (666,2006)(682,2023)(698,2040)
	 (715,2058)(733,2077)(752,2097)
	 (773,2118)(795,2141)(819,2165)
	 (844,2191)(870,2218)(897,2246)
	 (925,2273)(952,2301)(977,2326)
	 (1000,2350)(1020,2370)(1036,2386)
	 (1048,2398)(1056,2406)(1060,2410)(1062,2412)
 \path(837,12)(835,14)(830,19)
	 (822,27)(809,40)(792,58)
	 (770,80)(744,106)(715,135)
	 (684,167)(652,200)(619,234)
	 (586,267)(555,299)(525,330)
	 (496,360)(470,387)(445,413)
	 (422,438)(401,460)(382,482)
	 (363,502)(346,521)(330,539)
	 (314,557)(300,574)(283,595)
	 (266,615)(250,634)(235,654)
	 (220,675)(205,695)(191,715)
	 (177,736)(163,757)(151,778)
	 (138,800)(127,821)(116,842)
	 (106,863)(96,884)(88,904)
	 (80,924)(72,945)(66,965)
	 (60,984)(54,1004)(50,1024)
	 (45,1045)(41,1066)(37,1088)
	 (34,1111)(31,1134)(29,1158)
	 (27,1183)(25,1208)(24,1234)
	 (24,1261)(24,1287)(25,1313)
	 (26,1340)(28,1366)(31,1391)
	 (34,1416)(37,1440)(41,1463)
	 (46,1486)(51,1508)(56,1529)
	 (62,1549)(69,1570)(76,1590)
	 (84,1609)(92,1629)(101,1649)
	 (111,1669)(122,1690)(134,1710)
	 (146,1730)(159,1751)(172,1771)
	 (186,1791)(201,1811)(215,1831)
	 (230,1850)(245,1869)(260,1887)
	 (276,1905)(291,1923)(306,1940)
	 (321,1957)(337,1974)(351,1990)
	 (366,2006)(382,2023)(398,2040)
	 (415,2058)(433,2077)(452,2097)
	 (473,2118)(495,2141)(519,2165)
	 (544,2191)(570,2218)(597,2246)
	 (625,2273)(652,2301)(677,2326)
	 (700,2350)(720,2370)(736,2386)
	 (748,2398)(756,2406)(760,2410)(762,2412)
 \end{picture}
 }\end{minipage}}
 \caption{The intersection of two curved tubes}\label{fig:intersec}
 \end{figure}
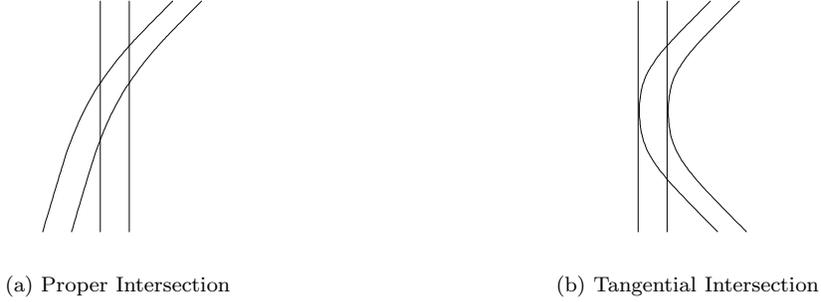

The proof of the ``trivial'' bound also
 reveals the reason for the non-degeneracy criterion \eqref{nondeg}, since if this does not hold, then the curves can essentially share a tangent, which makes the intersection of the tubes larger than the estimate given in Lemma~\ref{lem:intersec}.  This is illustrated in Figure~\ref{fig:intersec}.
As one would expect, this behaviour means that the $L^{\frac{n+1}{2}}$ estimate for the maximal function fails.
\begin{prop}If the non-degeneracy criterion \eqref{nondeg} fails, then \eqref{estimate} also fails for $p=\frac{n+1}{2}$.
\end{prop}
\begin{proof}
Suppose that a degenerate critical point occurs at $(t_0,z_0)$.  With $\psi (x_n,y):=\frac\partial{\partial x_n}\grad_y\tilde{\varphi}(x_n,y)$ as on page~\pageref{lem:intersec}, this means that $\det\frac{\partial}{\partial y}\psi(t,y)$ evaluated at $(t_0,z_0)$ is zero.  Let $U$ be the subspace
$$U=\left\{u\in\R^{n-1}:\frac{\partial}{\partial y}\psi(t_0,z_0)u=0\right\}$$ and let $r\geq1$ denote the dimension of this subspace.  As always, the curves have the parametrisation \eqref{xncurve}, and so for each $y$ we choose $\omega=\grad_y\tilde\varphi(t_0,y)$ to make all of the curves meet at the bad point.  
Consider those directions $y$ such that 
 $y-z_0\in U$. As before, the difference in tangent of the curves $\Gamma_y$ and $\Gamma_{z_0}$ at their intersection is $M(t_0)^{-1}[\psi(t_0,y)-\psi(t_0,z_0)]$.  By the definition of the derivative we then have
\begin{align*}
\psi(t_0,y)-\psi(t_0,z_0)&=\frac{\partial}{\partial y}\psi(t_0,z_0)(y-z_0)+O(|y-z_0|^2)\\
&=O(|y-z_0|^2)\quad\text{since $y-z_0\in U$,}
\end{align*}
which is at most $\delta$ provided that $|y-z_0|\lesssim\surd{\delta}$.
Pick a maximal $\surd\delta$-separated subset $\{y_j\}_{j=1}^M$ of these $y$.  Then $M\sim\surd\delta^{-r}$, and each tube $T_j:=T^\delta_y\big(\grad_y\tilde\varphi(t_0,y)\big)$ meets $T_0:=T^\delta_{z_0}\big(\grad_y\tilde\varphi(t_0,z_0)\big)$ at an angle of at most $\delta$.  So there is a  cylinder of radius $\delta$ and length some small constant $c$ which is included in all $M$ of the tubes.  Hence
\begin{align*}
\left\|\sum_{j=0}^M \1_{T_j}\right\|_{\frac{n+1}{n-1}}&\geq M (\delta^{n-1}c)^{\frac{n-1}{n+1}}\\
&\gtrsim\delta^{-\frac{n-1}{n+1}-\frac{r}{2(n+1)}}(\delta^{n-1}M)^{\frac{n}{n+1}}.
\end{align*}
Hence, by the covering lemma (Lemma~\ref{lem:covering}), we find that the $L^{\frac{n+1}{2}}\to L^{n+1}$ norm of $\kmf$ is at least $\delta^{-\frac{n-1}{n+1}-\frac{r}{2(n+1)}}$ which is greater than the estimate 
$\delta^{-\frac{n-1}{n+1}}$ obtained in the non-degenerate case.
\end{proof}
Rather curiously, the non-degeneracy criterion is not required for the set dimension problem.
\begin{prop}For any phase $\phi$ of the form \eqref{xnonlyphase} satisfying \eqref{rank} (but not necessarily \eqref{nondeg}), the corresponding curved Kakeya sets have Hausdorff and Minkowski dimension at least $\frac{n+1}{2}$.
\end{prop}
\begin{proof}This is intuitively clear, since a
 Kakeya set of degenerate curves includes a set of non-degenerate ones
 by simply removing slices around the ``bad'' heights.  Shifting and
 scaling part of what remains so that it lies in the region
 $x_n\in[-1,1]$ gives a set of curves that falls within the scope of
 Theorem~\ref{thm:trivial}.  So this subset, and hence the whole set, of
 the original curves has Minkowski and Hausdorff dimension at least
 $\frac{n+1}{2}$.  To prove this fact directly one merely needs to note that the conclusion of Lemma~\ref{lem:intersec} is true when $\delta=d$ whether the intersection is tangential or not, and then follow the proof of Theorem~\ref{thm:trivial}.\end{proof}
\section{Negative results for quadratic curves}
\label{sec:neg}
In the next three sections we look at the possibility of non-trivial results. Since we already know that these cannot hold of all classes of curves,  from here onwards restrict our attention to simpler ones.
Notice that in both of Bourgain's theorems the bad behaviour is caused by the presence of terms non-linear in $x$ in the phase function.  For this reason we now focus entirely on parabolic curves of the following form:
$$
\Gamma_y(\omega)=\left\{\begin{pmatrix}\omega-tAy-t^2By\\t\end{pmatrix}:t\in[-1,1]\right\}
$$
where $A$ and $B$ are $(n-1)\times(n-1)$ real symmetric matrices.  Kakeya questions about these curves arise from the phase
\begin{equation}
\label{easyphi}
\varphi(x,y):=y^\trns x'+x_n\frac12y^\trns\!A y+x_n^2\frac12y^\trns\!B y
\end{equation}
which is of the form \eqref{xnonlyphase}, and includes the Restriction problem as the special case $B=0$.  Meanwhile, the phase
\begin{align}
\varphi(x,y) &= \frac{1}{2x_n}y^\trns\left(A+x_nB\right)^{-1}y-\frac{1}{x_n}y^\trns\left(A+x_nB\right)^{-1}x'\label{nikphase}
\end{align}
is again of the form \eqref{xnonlyphase} and includes Bochner-Riesz as the case $B=0$, but it
gives rise to the same curves above but with $y$ and $\omega$ exchanged.  For this reason it is now convenient to call the Kakeya maximal function arising from \eqref{nikphase} a {\em Nikodym maximal function} and to denote it by
$$\nmf f(\omega):=\sup_{y\in\ball}
\frac1{|T_y(\omega)|}\int_{T_y(\omega)}|f(x)|\d x,$$ 
where we now fix the curves $\Gamma_y(\omega)$ and tubes $T_y(\omega)$ to be as above.  So $\Gamma_y(\omega)$ is always a parabola through the point $(\begin{smallmatrix}\omega\\0\end{smallmatrix})$ whose direction is governed by $y$.  Similarly, we define a {\em curved Nikodym set} to be one which includes a $\Gamma_y(\omega)$ for every $\omega\in\ball$.

It is easy to check that both phases above satisfy the first of H\"ormander's criteria \eqref{rank}; the second says we must assume that
\begin{equation}
\label{nondegAB}
\det(A+2x_nB)\neq 0
\end{equation}
throughout the support of the cutoff $a$, which we will take to be $[-1,1]$, with the obvious deletion of a neighbourhood of $x_n=0$ in the second case.
Also, by applying linear maps to $x$ and/or $y$ in the phase, we see that the
oscillatory integral problem is invariant under congruence of the matrices.  Since they are
symmetric, we may assume that one of them is diagonal, or even has
only $0$ and $\pm 1$ on the diagonal.  This is of limited help, but in
the special case where one of $A,B$ is positive-definite we are able
to simultaneously diagonalise (that is, using change of variable we can make both $A$ and $B$ diagonal).
  This will enable us to perform certain
computations that seem intractable in the general case.

At the level of the curves rather than the phases, we are free to multiply through by any invertible matrix,
This allows us, if it is convenient, to replace
$A$ by $I$ and $B$ by $C:=A^{-1}B$, so that the curves are now
\begin{equation}
\label{curves}
\Gamma_y(\omega)=\left\{\begin{pmatrix}\omega-ty-t^2Cy\\t\end{pmatrix}:t\in[-1,1]\right\}
\end{equation}
 Note however
that this matrix is not necessarily symmetric, nor is $B$
assumed to be invertible.  By a further transformation we may assume that $C$ is in rational canonical form.

Note also that if $B$ is a multiple of $A$ (so $C=\lambda I$ say\label{straighten}) then we can eliminate $C$ altogether using the diffeomorphism $x_n+\lambda x_n^2\mapsto x_n$.  So the curved case only arises if $C$ is not a multiple of $I$.
Moreover, in the case $C=0$, where the two phases above correspond to Restriction and Bochner-Riesz respectively, their corresponding maximal functions are related. This is because the transformation
$$(x',x_n)\mapsto\left(\frac{x'}{x_n},\frac{1}{x_n}\right)$$
(which was first used by Carbery \cite{carbery:paraboloid} in showing that Restriction implies Bochner-Riesz) maps straight lines to straight lines: specifically, the line 
 centred at $(\begin{smallmatrix}\omega\\0\end{smallmatrix})$ in direction $(\begin{smallmatrix}y\\1\end{smallmatrix})$ maps to the line with centre $(\begin{smallmatrix}y\\0\end{smallmatrix})$ and direction $(\begin{smallmatrix}\omega\\1\end{smallmatrix})$.  This is why the Kakeya and Nikodym problems for straight lines are equivalent.  Importantly, however, this transformation does {\em not} map parabolas to parabolas even if the roles of position and direction are exchanged.  So there is no reason to expect $\kmf$ and $\nmf$ to satisfy the same bounds for a given matrix $C$, and in fact we shall see later that they do not.

 Although very simple, these phases are general enough
to exhibit many kinds of behaviour.  Taking
$C=(\begin{smallmatrix}0&0\\1&0\end{smallmatrix})$ gives the ``worst case''
example of Theorem~\ref{thm:worst}.  More interestingly still, the
Generic Failure criterion \eqref{genericfailure} of Theorem~\ref{thm:genfail} has the simple form
\begin{equation*}
 C\mbox{ is not a multiple of }I.
\end{equation*}
Bourgain's proof of Theorem~\ref{thm:genfail} is considerably simpler for these special curves, works in higher dimensions, and in fact gives a better bound in dimension $3$, so we include the details here.
\begin{thm}
\label{thm:easygenfail}
Suppose that the characteristic polynomial of C divided by its minimum polynomial consists of irreducible factors each of multiplicity at most $k$, where $0\leq k\leq n-2$. Then 
$$\|\kmf\|_{p\to1}\gtrsim\delta^{\frac{1}{2n-3-k}(1-\frac1p)-\frac1p}$$ 
for all $p$. If $k=0$ and additionally $\tr\adj C=0$, then this is strengthened to $\delta^{\frac{1}{2n-2}(1-\frac1p)-\frac1p}$, while if $k=0$ and $\tr\adj C=0$ and $\det C=0$ it is strengthened further to $\delta^{-1/p}$. 
\end{thm}
This applies only to the Kakeya problem with parabolas and not the Nikodym version, although if $n=3$ then phases of the form \eqref{nikphase} are covered by Bourgain's generic failure result (Theorem~\ref{thm:genfail}) after making changes of variable and expanding the $1/x_n$ as a power series to obtain the standard form \eqref{generalphase}.

Combining this with the implication of Proposition~\ref{prop:oscillimplies} gives the following partial answers to Questions~\ref{que:max} and \ref{horconj}:
\begin{cor}
If  $\varphi$ is of the form \eqref{easyphi} with $B$
not a multiple of $A$, then the desired estimate $\|\kmf \|_{p\to q}\lesssim \delta^{1-n/p}$ for the corresponding Kakeya maximal function is false for all 
$p> n-\frac{n-k-2}{2n-k-4}$ even with $q=\infty$, and the estimate $\|T_N\|_{r\to s}\lesssim N^{-n/s}$ is false for all
$s<\frac{2n}{n-1}+\frac{2n-2k-2}{(2n-k-3)(n-1)(2n-3)}$
even with $r=1$.

If $k=0$ and additionally $\tr\adj C=0$, then these are strengthened to $p> n-\frac{n-1}{2n-3}$ and $s<\frac{2n}{n-1}+\frac{2n-2}{2(n-1)^2(2n-3)}$, while if $k=0$ and $\tr\adj C=0$ and $\det C=0$ these are strengthened further to $p> n-1$ and $s<\frac{2n}{n-1}+\frac{2}{(n-1)(2n-3)}$. 
\end{cor}
So we cannot achieve the optimal $p=n$, $s=\frac{2n}{n-1}$ unless $k=n-2$, which means that the minimum polynomial of $C$ is linear and so $C\parallel I$. Moreover, since a ``generic'' matrix has its characteristic and minimal polynomials equal, we can ``usually'' achieve no better than $p=n-\frac12$ and $s=\frac{2n}{n-1}+\frac{1}{(n-1)(2n-3)}$.

  Note that in three dimensions  we are dealing with $2\times2$ matrices and so we always have $k=0$ as long as $C\nparallel I$.
So we cannot exceed  the bound  $p=5/2=\frac{n+2}{2}$, which for straight lines is due to Wolff \cite{wolff:improvedbound}.  If additionally $\tr C=0$ we cannot exceed $p=7/3$, which was obtained for straight lines first by Bourgain \cite{bourgain:besicovitchtype} and then by Schlag \cite{schlag:geometricineq}.  Their results have since been improved for straight lines, but the above theorem for curves suggests that $5/2$ and $7/3$ are ``natural barriers'' in the problem.  The case $k=0$ and $\tr\adj C=0$ and $\det C=0$ in dimension three corresponds precisely to Bourgain's ``worst case'' example of Theorem~\ref{thm:worst}, and in fact there are analogues in higher dimensions as we shall see later.

The gain in three dimensions here compared with Theorem~\ref{thm:genfail} is because of the absence of higher order terms.  Their absence is also needed to make the proof work in dimensions $4$ and above, since Bourgain's method of handling the general case in \cite{bourgain:lp} uses that the order of the terms neglected is equal to the dimension $n$.
\begin{proof}[Proof of Theorem~\ref{thm:easygenfail}]It is enough to show that we can choose suitable
$\omega=\omega(y)$ to produce a set of curves that is too small.  We
will use a linear function: $\omega =Wy$. We claim that if we can make the determinant of the map $y\mapsto x':=Wy-ty-t^2Cy$ of
order $|t|^m$ for small $t$, then $\|\kmf\|_{p\to1}\gtrsim\delta^{1/m-1/p(1+1/n)}$ for all $p$.

Fix $t\in[-\delta^{1/m},\delta^{1/m}]$ so that the determinant is at most $\delta$.  Then if
 $y$ ranges over the ball in $\R^{n-1}$ of radius $1$, we find
 that $x'$ ranges over a set of measure at most $\delta$, and we are
 interested in the size of the $\delta$-neighbourhood of this.  Now
 since the eigenvalues of the map are bounded, no side of the set can
 be larger than $|y|<1$, but having all sides this large would exceed
 the maximum permitted volume.  So the worst case has $n-2$ sides of
 length $1$ and one thin side of length $\delta$ so that the volume
 does not exceed that permitted by the determinant.

Hence the largest possible neighbourhood is of measure $\delta$.  Now
allowing $x_n$ to vary over the interval
$[-\delta^{1/m},\delta^{1/m}]$ gives us that the union $E$ of these
tubes of length $\delta^{1/m}$ has measure at most $\delta^{1+1/m}$.  Observing that $\kmf\1_E(y)\geq\delta^{1/m}$ for all $y\in\ball$ proves the claim.

So we must consider when the above condition on the determinant is
satisfied for some $m\geq n$.   Clearly if $C\parallel I$ then it cannot be, since the determinant is just the characteristic polynomial of $W$ evaluated at $(t+\lambda t^2)$, but in all other cases we can simply write down a suitable $W$.  By a change of variable we may assume that $C$ is in rational canonical form; that is, $C =C_{p_0}\oplus\dots \oplus C_{p_k}$ where $k$ is as above, each $C_{p_i}$ is the companion matrix of the polynomial $p_i$,  $p_k$ is the minimum polynomial of $C$, and $p_i$ divides $p_{i+1}$ for $i=0,\dots,k-1$. 

We show that for an $l\times l$ companion matrix 
$$C=\begin{pmatrix}c_1 &1 &0 &0&\dots &0\\
c_2 &0&1&0&\dots &0\\
c_3 &0&0&1&\dots&0\\
\vdots&\vdots&\vdots &\ddots&\ddots &\vdots\\
c_{l-1}&0&0&\dots &0&1\\
c_{l}&0&0&\dots&0&0
\end{pmatrix}$$
we can achieve $\det(W-tI-t^2C)\lesssim t^{2l-1}$.  Choose $W$ to be zero except in the first column, whose elements are as follows:
\begin{align*}
w_{1,1}&=0& w_{2,1}&=-1&w_{i,1}&=c_{i-2}\text{ for }i\geq3.
\end{align*}
If we expand $\det(W-2tI-2t^2C)$ down the first column we obtain
\begin{align*}&(-1)^l\left[(c_1t^2+t)t^{l-1}-(c_2t^2+1)t^l+\sum_{i=3}^l(-1){i-1}(c_it^2-c_{i-2})t^{l+i-2}\right]\\
=\;&(-1)^l\left[c_1t^{l+1}-c_2t^{l+2}+\sum_{i=3}^l(-1)^{i-1}c_it^{l+i}-\sum_{i=1}^{l-2}(-1)^{i-1}c_it^{l+i}\right]\\
=\;&c_{l-1}t^{2l-1}-c_lt^{2l}
\end{align*}
This proves the result since there are $k+1$  blocks each of order $l_i\times l_i$ with $\sum_i l_i=n-1$, and so we take $m=2(n-1)-k-1$. If $k=0$ then we have just one block of order $(n-1)\times (n-1)$. The conditions $\tr\adj C=0$ and $\det C=0$ correspond to $c_{n-2}=0$ and $c_{n-1}=0$, which allow us to take $m=2(n-1)$ or $m=\infty$ respectively, giving the improvements stated.
\end{proof}
\section{Geometric Methods}
\label{sec:geometric}
In this section we will prove a result for the curved Kakeya maximal function of a particular class of curves which implies that the corresponding sets have dimension at least $\frac{n+2}{2}$.  In the straight line case this result is due to Wolff \cite{wolff:improvedbound} and uses geometric techniques, however, here we shall adapt a more recent proof due to Katz \cite{katz:socialdensity}.  In both, the main geometric object is the {\em hairbrush}, a configuration of tubes which all pass through some central fixed one.  Wolff's idea was that such configurations can be handled by grouping the tubes into planes all containing the central tube, and then by applying the known dimension $2$ result for the Kakeya maximal function in each plane.  Curves, however, cannot easily be grouped in this way, which is why we turned to Katz's work.  His proof seems more
 elementary, in that it isolates the geometry showing that the main
 fact is that a triangle lies in a plane, and the remainder of the
 argument is a simple (but clever) splitting up of the linearised
 maximal function into bounded pieces.

 Our result is the following:
 \begin{thm}
 \label{thm:katz}
 The Nikodym maximal function $\nmf$ satisfies the bound 
 $$\|\nmf\|_{n\to n}\lesssim\delta^{-\frac{n-2}{2n}}$$ provided that the curves under consideration are parabolas of the form \eqref{curves} with $C^2=0$.
 \end{thm}
By Proposition~\ref{prop:maxtohaus} this implies that the Nikodym sets of these curves have Hausdorff and Minkowski dimension at least $\frac{n+2}{2}$.

The condition on the matrix $C$ arises naturally in the
 proof as we shall see in Proposition~\ref{prop:property}.  This class of curves
 seems to be particularly amenable to the proof methods that have
 been used in the straight line case, since further results for these
 curves will be obtained by the arithmetic methods in Section~\ref{sec:arith}.
 However, the class is not equivalent to straight lines since, as we will see in Section~\ref{sec:discuss}, the Kakeya conjecture fails completely for all these curves.  Unfortunately it seems difficult to give the criterion $C^2=0$ any geometric interpretation.

 We shall actually prove Theorem~\ref{thm:katz} for the linearised
 version of the Nikodym maximal function $\lnmf$---recall its
 definition from page~\pageref{def:LK}: We have divided $\R^{n-1}$ into
 $\delta$-cubes $Q_j$ where $j$ runs over $\ball\cap\delta\Z^{n-1}$.
 To each index $j$ we have an associated curved tube
 $T_j=T_{y_j}(\omega_j)$ where $\omega_j\in Q_j$ and $y_j$ is
 arbitrary.  Then
 $$\lnmf f(\omega)=\sum_j\1_{Q_j}(\omega)\frac{1}{|T_j|}\int_{T_j}f(x)\d x.$$
Of course, we must seek bounds that are independent of the choice of
the tubes.  We shall also need to define related functions where the
index set is specified:
 $$\Lcal\Ncal_\Abb f(\omega)=\sum_{j\in\Abb}\1_{Q_j}(\omega)\frac{1}{|T_j|}\int_{T_j}f(x)\d x.$$

 As in Wolff's approach, the main geometric object considered is the {\em hairbrush}:
 \begin{defn}
 Let $\Abb$ be a finite set of indices $j\in\delta\Z^{n-1}$.  A {\em hairbrush} is a set $\H\subseteq \Abb$ such that there exists some curved $1\times\delta$ tube $T$ that intersects all $T_i$ with $i\in\H$.
 \end{defn}
\noindent Note that the central tube $T$ can be any curved tube of the family, not necessarily one of those associated to some $j$.

 Much of the geometry of the situation is encoded in the behaviour of
 these hairbrushes, in the form of the following lemma:
 \begin{lem}[Hairbrush Lemma]\label{lem:hairbrush} If the curves are
 parabolas with $C^2=0$, then for all hairbrushes $\H$ we have
 $\|\Lcal\Ncal_\H\|_{n\to n}\leq C(\log1/\delta)^\alpha$ for some constant $\alpha$  \end{lem}
 The proof of this will involve surfaces, and will show why we are
 able to handle only a restricted class of curves.  But given the
 lemma, we can prove the theorem just as in the straight line case, by
 splitting up the operator into many sums.
\begin{proof}[Proof of Theorem~\ref{thm:katz}] 
 As usual it is enough to prove a weak type estimate.  By the covering lemma (Lemma~\ref{lem:covering}) the theorem is true if and only if 
 \begin{align*}
 &&\left\|\sum_{j\in\Abb}\1_{T_j}\right\|_{n'}&\lesssim \delta^{-\frac{n-2}{2n}}(\delta^{n-1}\#\Abb)^{1/n'}\\
 &\iff & \int_{\R^n} \Bigl(\sum_{j\in\Abb}\1_{T_j}(x)\Bigr)^{\frac{n}{n-1}}\d x &\lesssim \delta^{-\frac{n-2}{2(n-1)}}(\delta^{n-1}\#\Abb)\\
 &\iff & \sum_{j\in\Abb}\frac{1}{|T_j|}\int_{T_j}\Bigl(\sum_{i\in\Abb}\1_{T_i}(x)\Bigr)^{\frac{1}{n-1}}\d x&\lesssim\delta^{-\frac{n-2}{2(n-1)}}\#\Abb
 \end{align*}
 Denote the quantity appearing in the first sum by $M_\Abb$, that is
 $$M_\Abb(j):=\frac{1}{|T_j|}\int_{T_j}\Bigl(\sum_{i\in\Abb}\1_{T_i}(x)\Bigr)^{\frac{1}{n-1}}\d x.$$
We would like to subdivide this quantity into dyadic scales, by
considering those $i$ that are at distance between $2^{-k}$ and
$2^{-(k+1)}$ from $j$.  Note that by elementary properties of
sequences of positive reals, the sum over $k$ can then be pulled out
of the integral.  But what remains then depends only on pairs $i,j$
with $|i-j|\sim 2^{-k}$.  Let $P$ be a cube in $\R^{n-1}$ of side
$10\times2^{-k}$.  Note that this cube is larger than the cubes $Q_j$
since $\delta<2^{-k}$.  It then suffices, for every choice of $P$, to
obtain the estimate
\begin{align*}\sum_{j\in\Abb\cap P}M_k(j)&:=\sum_{j\in\Abb\cap P}\frac{1}{|T_j|}\int_{T_j}\Big(\sum_{\substack{i\in\Abb\cap P\\|i-j|\sim2^{-k}}}\1_{T_i}(x)\Big)^{\frac{1}{n-1}}\d x
\\& \lesssim\delta^{-\frac{n-2}{2(n-1)}}\#(\Abb\cap P)\end{align*}
 and then sum over $P$ and $k$.  Both sums have only logarithmically many terms.  

 The next stage is to find as many large hairbrushes in $\Abb\cap P$
 as possible, where {\em large} means of cardinality at least $N$, to
 be chosen later.  So, if there exists some curved tube $T$ (of the
 form $T_y(\omega)$ but not necessarily one of the $T_j$) such that
 there are at least $N$ elements $i\in\Abb\cap P$ with $T\cap
 T_i\neq\emptyset$, then call these elements $\H_1$.  Then look for
 another large hairbrush in the remaining elements $\Abb\cap
 P\setminus\H_1$.  Eventually there are no more hairbrushes, so
 call the remaining bad elements $\Bbb$.  This constructs hairbrushes
 $\H_1,\dots,\H_m$ each of cardinality at least $N$, and a bad set
 $\Bbb:=\Abb\setminus(\H_1\cup\dots\cup\H_m)$.  Let
 $\H:=\H_1\cup\dots\cup\H_m$.

 Since the hairbrushes are disjoint sets of indices (although the
 tubes they correspond to may well not be), and $\Abb\cap P$ has at
 most $2^{-k(n-1)}\delta^{-(n-1)}$ elements, it follows that $m\leq
 2^{-k(n-1)}\delta^{-(n-1)}/N$.

 Now split the sum into four pieces
 $$\sum_{j\in\Abb\cap P}M_k(j)\leq \sum_{j\in\H}M_{k,\H}(j)+\sum_{j\in\Bbb}M_{k,\H}(j)+\sum_{j\in\H}M_{k,\Bbb}(j)+\sum_{j\in\Bbb}M_{k,\Bbb}(j)$$
 where
 \begin{align*}
 M_{k,\H}(j)&:=\frac{1}{|T_j|}\int_{T_j}\Big(\sum_{\substack{i\in\H\\|i-j|\sim2^{-k}}}\1_{T_i}(x)\Big)^{\frac{1}{n-1}}\d x\\ M_{k,\Bbb}(j)&:=\frac{1}{|T_j|}\int_{T_j}\Big(\sum_{\substack{i\in\Bbb\\|i-j|\sim2^{-k}}}\1_{T_i}(x)\Big)^{\frac{1}{n-1}}\d x.
 \end{align*}

 The first sum is estimated using the hairbrush lemma.  For
 \begin{align*}
 \|\Lcal\Ncal_\H f\|^n_n & = \int\Bigl(\sum_i \Lcal\Ncal_{\H_i}f(\omega)\Bigr)^n\d \omega \\
 & = \sum_i\int\bigl( \Lcal\Ncal_{\H_i}f(\omega)\bigr)^n\d \omega &\mbox{since each $\omega$ gives only one non-zero term}\\
 & \leq m\|\Lcal\Ncal_{\H_i}\|^n_{n\to n}\|f\|^n_n
 \end{align*}
 showing that $\|\Lcal\Ncal_\H\|_{n\to n}\leq Cm^{1/n}(\log1/\delta)^\alpha$ by Lemma~\ref{lem:hairbrush}.  Then by the covering lemma we obtain
 $$\sum_{j\in\H}M_{k,\H}(j)\leq (Cm^{1/n}(\log1/\delta)^\alpha)^{n'}\#\H.$$

 For the second sum
 \begin{alignat*}{2}
 \sum_{j\in\Bbb}M_{k,\H}(j) &=\sum_{j\in\Bbb}\frac{1}{|T_j|}\int_{T_j}\Big(\sum_{\substack{i\in\H\\|i-j|\sim2^{-k}}}\1_{T_i}(x)\Big)^{\frac{1}{n-1}}\d x\\
 &\leq \sum_{j\in\Bbb}\Bigg(\frac{1}{|T_j|}\int_{T_j}\sum_{\substack{i\in\H\\|i-j|\sim2^{-k}}}\1_{T_i}(x)\d x\Bigg)^{\frac{1}{n-1}} &\mbox{by Jensen}\\
 & \leq (\#\Bbb)^{1-\frac{1}{n-1}}\Bigg(\sum_{j\in\Bbb}\frac{1}{|T_j|}\int_{T_j}\sum_{\substack{i\in\H\\|i-j|\sim2^{-k}}}\1_{T_i}(x)\d x\Bigg)^{\frac{1}{n-1}} &\mbox{by H\"older}\\
 & \leq (\#\Bbb)^{1-\frac{1}{n-1}}\Bigg(\sum_{i\in\H}\sum_{\substack{j\in\Bbb\\|i-j|\sim2^{-k}}}\frac{1}{|T_j|}|T_i\cap T_j|\Bigg)^{\frac{1}{n-1}} &\mbox{by swapping sums}\\
 & \lesssim (\#\Bbb)^{1-\frac{1}{n-1}}\Bigg(\sum_{i\in\H}\sum_{\substack{j\in\Bbb\\T_i\cap T_j\neq\emptyset}}\frac{1}{\delta^{n-1}}\delta^n2^k\Bigg)^{\frac{1}{n-1}} &\mbox{by Lemma~\ref{lem:intersec}}\\
 & \leq (\#\Bbb)^{1-\frac{1}{n-1}}\left(\sum_{i\in\H}N\delta 2^k\right)^{\frac{1}{n-1}}&\mbox{\hspace{-4cm}since no large hairbrushes in $\Bbb$}\\
 & =(\#\Bbb)^{1-\frac{1}{n-1}}\left(\#\H N\delta 2^k\right)^{\frac{1}{n-1}} \\
 & \leq (\#\Abb\cap P)\big( N\delta 2^k\big)^{\frac{1}{n-1}}. 
 \end{alignat*}

 The third and fourth sums can be tackled together, since for {\em all} $j\in\Abb\cap P$ we have
 \begin{alignat*}{2}
 M_{k,\Bbb}(j)&:=\frac{1}{|T_j|}\int_{T_j}\Big(\sum_{\substack{i\in\Bbb\\|i-j|\sim2^{-k}}}\1_{T_i}(x)\Big)^{\frac{1}{n-1}}\d x\\
 & \leq \Bigg(\frac{1}{|T_j|}\int_{T_j}\sum_{\substack{i\in\Bbb\\|i-j|\sim2^{-k}}}\1_{T_i}(x)\d x\Bigg)^{\frac{1}{n-1}}&\mbox{by Jensen}\\
 & = \Big(\frac{1}{|T_j|}\sum_{\substack{i\in\Bbb\\|i-j|\sim2^{-k}}}|T_i\cap T_j|\Big)^{\frac{1}{n-1}}&\mbox{}\\
 &\lesssim  \Big(\frac{1}{\delta^{n-1}}\sum_{\substack{i\in\Bbb\\T_i\cap T_j\neq\emptyset}}\delta^n2^k\Big)^{\frac{1}{n-1}}&\mbox{by Lemma~\ref{lem:intersec}}\\
 &\leq  \big(N\delta 2^k\big)^{\frac{1}{n-1}}&\mbox{since no large hairbrushes in $\Bbb$}.
 \end{alignat*}

 So the last three sums all give an estimate of $\#(\Abb\cap
 P)(N\delta 2^k)^{\frac{1}{n-1}}$, while the first, after putting in
 the upper bound for $m$, gives $\#\H
 (2^{-k(n-1)}\delta^{-(n-1)}/N)^{\frac{1}{n-1}}$.  We optimally
 choose $N=2^{-k\frac{n}{2}}\delta^{-\frac n2}$ and add the four
 pieces to obtain
 $$\sum_{j\in\Abb\cap P}M_k(j)\lesssim 2^{-k\frac{n-2}{2(n-1)}}\delta^{-\frac{n-2}{2(n-1)}}\#(\Abb\cap P)$$
 which gives the result after summing over all $P$ of side $2^{-k}$ and all $k$.\end{proof}
To prove the hairbrush lemma for curves, we need an analogue of the following fact about straight lines: 
\begin{quote}\label{fact}Two intersecting straight lines determine a plane, and thus a third line intersecting these two is fixed up to one parameter, i.e.~its direction must lie parallel to the plane, or equivalently the point where it meets the base plane $x_n=0$ must lie along a fixed line.\end{quote}  So we must now study the locus of all curves meeting two given ones.

 By the linearity of \eqref{curves} in $y$ and $\omega$, we may assume
 that one of the given curves is $\Gamma_0(0)$. Let the other be
 $\Gamma_{y_0}(\omega_0)$ and assume that they meet at height
 $t_0$. The surface is the locus of those curves $\Gamma_y(\omega)$
 that meet the first at $s$ and the second at $u$.  Note that none of
 these three heights are equal, since the curves are never horizontal,
 and we must exclude the possibility of $\Gamma_y(\omega)$ meeting the
 two given curves at their common point, since this would allow every
 curve to belong to the locus. This is made clearer by the following
 picture:
 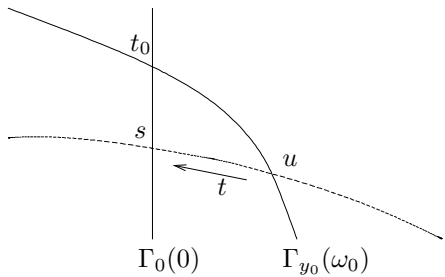
\begin{figure}[!hbtp]
 \centering
 \setlength{\unitlength}{0.0005in}
 {\renewcommand{\dashlinestretch}{30}
 \begin{picture}(4524,2781)(0,-10)
 \path(2490,968)(1740,1118)
 \path(1863.553,1123.883)(1740.000,1118.000)(1851.786,1065.049)
 \path(1512,2754)(1512,354)
 \dashline{60.000}(4512,354)(4511,354)(4510,355)
	 (4507,357)(4503,359)(4496,362)
	 (4488,366)(4477,372)(4463,379)
	 (4447,387)(4428,396)(4407,406)
	 (4383,418)(4357,431)(4329,445)
	 (4298,459)(4266,475)(4231,491)
	 (4196,508)(4158,525)(4119,543)
	 (4079,561)(4037,580)(3994,598)
	 (3949,618)(3903,637)(3855,657)
	 (3806,677)(3754,698)(3700,719)
	 (3644,740)(3585,762)(3524,785)
	 (3459,808)(3391,831)(3321,855)
	 (3247,880)(3171,905)(3092,929)
	 (3012,954)(2935,977)(2859,999)
	 (2786,1020)(2715,1040)(2648,1058)
	 (2586,1075)(2528,1090)(2475,1104)
	 (2427,1117)(2384,1128)(2345,1137)
	 (2311,1146)(2280,1153)(2253,1159)
	 (2230,1165)(2209,1169)(2190,1173)
	 (2174,1177)(2159,1180)(2144,1183)
	 (2131,1185)(2117,1188)(2103,1190)
	 (2088,1193)(2072,1196)(2054,1199)
	 (2033,1202)(2010,1206)(1984,1211)
	 (1955,1216)(1921,1222)(1884,1229)
	 (1842,1236)(1796,1244)(1745,1253)
	 (1690,1263)(1630,1273)(1567,1284)
	 (1500,1295)(1430,1307)(1359,1318)
	 (1287,1329)(1200,1342)(1116,1353)
	 (1035,1364)(959,1373)(887,1381)
	 (820,1387)(757,1393)(697,1398)
	 (642,1402)(589,1405)(540,1408)
	 (492,1410)(448,1412)(405,1413)
	 (364,1413)(326,1413)(289,1413)
	 (253,1413)(220,1413)(189,1412)
	 (160,1411)(134,1410)(110,1409)
	 (89,1408)(70,1407)(54,1406)
	 (41,1406)(31,1405)(23,1405)
	 (18,1404)(15,1404)(13,1404)(12,1404)
 \path(3012,354)(3011,356)(3010,361)
	 (3007,369)(3002,382)(2996,401)
	 (2987,425)(2977,453)(2965,487)
	 (2951,524)(2937,564)(2921,607)
	 (2905,650)(2889,694)(2873,738)
	 (2856,780)(2841,821)(2825,860)
	 (2810,897)(2796,933)(2781,966)
	 (2768,998)(2754,1027)(2741,1056)
	 (2728,1083)(2715,1108)(2702,1133)
	 (2689,1157)(2676,1181)(2662,1204)
	 (2647,1228)(2632,1253)(2616,1277)
	 (2599,1301)(2582,1325)(2565,1349)
	 (2546,1373)(2527,1398)(2507,1422)
	 (2486,1447)(2465,1471)(2443,1495)
	 (2421,1520)(2398,1544)(2375,1567)
	 (2351,1591)(2327,1613)(2303,1636)
	 (2279,1657)(2255,1679)(2231,1699)
	 (2207,1719)(2183,1739)(2159,1758)
	 (2135,1776)(2111,1794)(2086,1812)
	 (2062,1829)(2039,1845)(2015,1861)
	 (1991,1877)(1967,1893)(1941,1909)
	 (1915,1925)(1889,1941)(1861,1958)
	 (1833,1974)(1803,1991)(1774,2008)
	 (1743,2025)(1712,2042)(1680,2059)
	 (1648,2076)(1616,2092)(1583,2109)
	 (1550,2125)(1517,2142)(1484,2158)
	 (1451,2174)(1419,2189)(1386,2204)
	 (1354,2219)(1321,2234)(1289,2248)
	 (1258,2262)(1226,2276)(1194,2290)
	 (1162,2304)(1134,2316)(1105,2328)
	 (1076,2341)(1046,2353)(1015,2366)
	 (983,2379)(950,2392)(916,2406)
	 (880,2420)(842,2435)(803,2451)
	 (762,2467)(718,2484)(673,2501)
	 (626,2520)(577,2538)(526,2558)
	 (475,2578)(423,2598)(371,2617)
	 (320,2637)(271,2656)(225,2673)
	 (181,2690)(142,2705)(108,2717)
	 (80,2728)(57,2737)(39,2744)
	 (26,2749)(18,2752)(14,2753)(12,2754)
 \put(2187,804){$t$}
 \put(1320,1404){$s$}
 \put(1362,54){$\Gamma_0(0)$}
 \put(2862,54){$\Gamma_{y_0}(\omega_0)$}
 \put(1280,2310){$t_0$}
 \put(2862,1104){$u$}
 \end{picture}
 }
\caption{Notation for a curved triangle}\label{fig:curvedtriangle}
 \end{figure} 
 \\Now we have the following equations:
 \begin{align}
 0& =  \omega_0-t_0y_0-t_0^2Cy_0\label{first}\\
 0 & =  \omega-sy-s^2Cy\label{second}\\
  \omega_0-uy_0-u^2Cy_0 & =  \omega-uy-u^2Cy\label{third}
 \end{align}
 Subtracting \eqref{second} from \eqref{third} we find that
 $$y=\frac{1}{s-u}\big(I+(s+u)C\big)^{-1}\big(\omega_0-uy_0-u^2Cy_0\big)$$
 which is well defined because $s\neq u$ and by \eqref{nondegAB}.
 Substitute into \eqref{second} to find $\omega$:
 $$\omega = \frac{s}{s-u}\big(I+sC\big)\big(I+(s+u)C\big)^{-1}\big(\omega_0-uy_0-u^2Cy_0\big).$$
 \eqref{first} has not been used yet, so we use it to eliminate $\omega_0$:
 $$\omega = \frac{s(t_0-u)}{s-u}\big(I+sC\big)\big(I+(s+u)C\big)^{-1}\big(I+(t_0+u)C\big)y_0.$$
 Finally substitute this $y$ and $\omega$ into \eqref{curves} to obtain
 \begin{equation}
 \label{fatsurf}
 \begin{pmatrix}\frac{(s-t)(t_0-u)}{s-u}\big(I+(s+t)C\big)\big(I+(s+u)C\big)^{-1}\big(I+(t_0+u)C\big)y_0\\t\end{pmatrix}
 \end{equation}
 as the parametrisation of the locus we are interested in.  Note that if $C=0$ then this reduces to the plane
 $(\begin{smallmatrix}ry_0\\t\end{smallmatrix})$ as expected.  In general however, there are three parameters $(u,s,t)$ and so the locus is not a plane nor even a surface but rather some fat object.  What we need to know is whether a curve belonging to this locus has its direction and/or its base point fixed up to one parameter.
The
 following proposition determines when this is so.
\begin{prop}
\label{prop:property}
Suppose that the curve $\Gamma_y(\omega)$ is included in the locus \eqref{fatsurf}.  Then
\begin{enumerate}
\item The point $\omega$ must belong to a family described by only one parameter if and only if $C\parallel I$ or $C^2=0$, in which cases 
\begin{equation}\label{nonfatomega}\omega=r(I+t_0C)y_0\end{equation}
 for some parameter $r$.
\item The direction $y$ must belong to a family described by only one parameter if and only if we have the straight line case $C\parallel I$.
\end{enumerate}
\end{prop}
\begin{proof}
From 
$$\omega-t(I+tC)y=\frac{(s-t)(t_0-u)}{s-u}\big(I+(s+t)C\big)\big(I+(s+u)C\big)^{-1}\big(I+(t_0+u)C\big)y_0$$
we obtain 
\begin{align*}
\omega &= \frac{s(t_0-u)}{s-u}(I+sC)\big(I+(s+u)C\big)^{-1}\big(I+(t_0+u)C\big)y_0\\
y &= \frac{t_0-u}{s-u}\big(I+(s+u)C\big)^{-1}\big(I+(t_0+u)C\big)y_0
\end{align*}
where for a given curve, $s$ and $u$ will be fixed.  However, we are considering the sets of all such $y$ and $\omega$, so we allow $s$ and $u$ to vary.  We also require the property for all $y_0$ and $t_0$.
\begin{proofenum}
\item To show that the locus of all $\omega$ is one-dimensional we require that the derivatives of $\omega$ with respect to $s$ and $u$ are always parallel.  These are
 \begin{multline*}
 \frac{t_0-u}{s-u}
 \bigg[s\big(C-(I+sC)(I+(s+u)C)^{-1}C\big)\\
 -\frac{u}{s-u}\big(I+sC\big)
 \bigg]\big(I+(s+u)C\big)^{-1}\big(I+(t_0+u)C\big)y_0
 \end{multline*}
 and
 \begin{multline*}
 \frac{s}{s-u}\big(I+sC\big)\big(I+(s+u)C\big)^{-1}\bigg[\frac{t_0-s}{s-u}\big(I+(t_0+u)C\big)+\\(t_0-u)\big(C-C(I+(s+u)C)^{-1}(I+(t_0+u)C)\big)\bigg]y_0.
 \end{multline*}
We need this for all $y_0$, so
 that in fact the matrices themselves must be ``parallel'', by which we
 mean that one is a scalar multiple of the other.  Next we may rewrite
 the above, but ignore the initial (scalar) function of $(s,u)$ and
 multiply on the left by $\big(I+(s+u)C\big)\big(I+sC\big)^{-1}$ and on the right by
 $\big(I+(t_0+u)C\big)^{-1}\big(I+(s+u)C\big)$.  We thus require the following two
 expressions to be parallel:
 \begin{gather*}
 \frac{1}{s-u}\big(I+(s+u)C\big)-sC\big(I+sC\big)^{-1}C\\
 \frac{1}{s-u}\big(I+(s+u)C\big)-(t_0-u)C\big(I+(t_0+u)C\big)^{-1}C
 \end{gather*}
Setting $s=0$ and $u=-t_0$ we find that $I-t_0C$, which is invertible, is parallel to $I-t_0C-2t_0C^2$.  This implies that either $C\parallel I$ or $C^2=0$.
\item For $y$, the two derivatives are
\begin{gather*}
\frac{t_0-u}{s-u}\big(I+(s+u)C\big)^{-1}\Big[-C\big(I+(s+u)C\big)^{-1}-\frac{1}{s-u}I\Big]\big(I+(t_0+u)C\big)\\
\begin{split}
\frac{1}{s-u}\big(I+(s+u)C\big)^{-1}\bigg[\frac{t_0-s}{s-u}\big(I+(t_0&+u)C\big)\\&+(t_0-u)C\big(I-\big(I+(s+u)C\big)^{-1}\big(I+(t_0+u)C\big)\big)\bigg]\end{split}.
\end{gather*}
Ignoring scalar functions and multiplying by invertible matrices on the right and left, we thus require the following two expressions to be parallel:
\begin{gather*}
I+2sC\\
I+(s+u)C-(t_0-u)(s-u)C\big(I+(t_0+u)C\big)^{-1}C
\end{gather*}
Setting $u=-t_0$ and subtracting gives
$$I+2sC\parallel  (s+t_0)C\left[I+2t_0C\right]$$
and since $I+2sC$ and $I+2t_0C$ are invertible by \eqref{nondegAB}, we can deduce that $C$ is a (possibly zero) multiple of $I$.
\end{proofenum}
In order to convince ourselves, we check that if $C^2=0$, then $y$ is given by $\frac{2(t_0-u)}{s-u}\big(I+(t_0-s)C\big)y_0$, which does have two parameters unless $C\parallel I$.
\end{proof}
This result clearly shows that for parabolas, the Kakeya and Nikodym versions of the problem are not the same at all.  Indeed, the ``worst case'' example of Bourgain had $C^2=0$, and we already know that no non-trivial Kakeya estimate can hold for this.

 We are now ready to prove the Hairbrush Lemma, and hence complete the proof of Theorem~\ref{thm:katz}.\begin{proof}[Proof of Lemma~\ref{lem:hairbrush}] 
 We have a set $\H$ of indices which forms a hairbrush with central
 tube $T$.  By linearity assume that the central tube is
 $T_0(0)$.  Denote the other tubes by $T_j=T_{y_j}(\omega_j)$, where
 $\omega_j\in Q_j$ and so $\omega_j\approx j$.  We partition the
 set $\H$ in several ways.  First, let $H_k$ be the set of all those
 indices whose tubes meet $T$ at ``angle'' $2^{-k}$; that is,
 $$H_k:=\{i\in\H: |y_i|\sim 2^{-k}\}.$$ 
For a fixed $\omega$, there
can be only one $k$ such that $\Lcal\Ncal_{H_k}(\omega)\neq0$, so it
is enough to prove $\|\Lcal\Ncal_{H_k}\|\leq C(\log1/\delta)^\alpha$
since there are only logarithmically many $k$.  Then by the arguments
used previously, this bound is true if and only if
 $$\sum_{i\in H_k}M_{H_k}(i)\leq C(\log1/\delta)^\alpha\#H_k.$$
 For fixed $j\in H_k$ split up $H_k$ into further sets $H_{j,k,l,m}$ as follows:
 $$H_{j,k,l,m}:=\left\{i\in H_k:|y_i-y_j|\sim2^{-l}\text{ and }\dist(T_i\cap T_j,T_j\cap T)\sim \delta2^{l+m}\right\}.$$
 Note that this set is empty unless $l\geq k-2$.  Now it is enough to show that $M_{H_{j,k,l,m}}(j)\leq C(\log1/\delta)^\alpha$, because 
 $$M_{H_k}(i)\leq\sum_{j\in H_k}\sum_l\sum_m M_{H_{j,k,l,m}}(i)$$
 and there are only logarithmically many $l$ and $m$ and the sum over $j$ introduces a factor $\#H_k$.

 Next comes the geometric part of the argument, which is a quantitative version of the fact explained on page~\pageref{fact}.  We need to show that
 $\#\H_{j,k,l,m}$ is not too big, which means that given the central
 tube $T$ and another fixed tube $T_j:j\in H_k$ there are few other tubes $T_i$ meeting these with all the correct ``angles''
 and distances.  In the straight line case this follows from simple
 consideration of similar triangles as in Figure~\ref{straightkatz}.
 \begin{figure}[!hbtp]
 \centering
 \setlength{\unitlength}{0.0006in}
 {\renewcommand{\dashlinestretch}{30}
 \begin{picture}(3396,3059)(0,-10)
 \dashline{60.000}(150,922)(2250,22)
 \path(150,322)(2850,3022)
 \path(1301.826,1292.225)(1200.000,1222.000)(1322.893,1236.045)
 \path(1200,1222)(2400,1672)
 \path(2298.174,1601.775)(2400.000,1672.000)(2277.107,1657.955)
 \path(150,922)(3150,2122) 
 \thicklines
 \path(2250,3022)(2250,22)
 \put(2290,2650){$2^{-k}$}
 \put(2290,1850){$2^{-k}$}
 \put(1500,1522){$2^{-l}$}
 \put(1650,1222){$\delta2^{l+m}$}
 \put(0,922){$j$}
 \put(300,172){$i$}
 \put(2400,22){$T$}
 \put(3150,1822){$T_j$}
 \end{picture}
 }
 \caption[{\em Illustration of a lemma on triangles}]{In the straight line case, by similar triangles we have $|i-j|\lesssim 2^{-l}$ and $\dist(i,\mathrm{line})\lesssim 2^{-(l+m)}$.}\label{straightkatz}
 \end{figure}
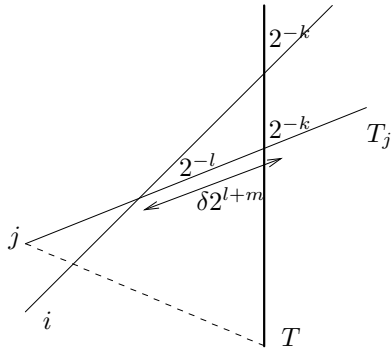  \\
In the curved case, the dotted line in the picture is instead the
line \eqref{nonfatomega}, which is the intersection of the base plane ${x_n=0}$ with the surface \eqref{fatsurf}
determined by $T$ and $T_j$.  Since we cannot appeal to similar
triangles with curves, we state and prove our claim more formally:
\begin{claim}\label{claim}
 Let $C$ satisfy $C^2=0$.  Suppose that we are given three curved tubes $T=T_0(0)$, $T_j=T_{y_j}(\omega_j)$ and $T_i=T_{y_i}(\omega_i)$ with $|y_j|,\,|y_i|
 \in(\frac{1}{2^{k+1}},\frac{1}{2^{k}})$ and $|y_j-y_i|\in(\frac{1}{2^{l+1}},\frac{1}{2^{l}})$.  Here $l\geq k-2$ and all the powers of $2$ that occur are greater than $\delta$.  Suppose that $T_j$ meets the axis at height $t_j$, $T_i$ meets it at $t_i$, and they meet each other at $s$, where $\delta2^{l+m}\leq |s-t_j|\leq \delta2^{l+m+1}$.  Then $|\omega_j-\omega_i|\leq\frac{1}{2^l}$, and $\omega_i$ is at distance at most $\frac{1}{2^{l+m}}$ from the line \eqref{nonfatomega}.
 \end{claim}

We have the following equations
 \begin{align*}
 \omega_j &= t_j(I+t_jC)y_j\\
 \omega_i &=t_i(I+t_iC)y_i+\eps\\
 \omega_i-s(I+sC)y_i & = \omega_j-s(I+sC)y_j+\eta
 \end{align*}
 where $\eps$ and $\eta$ are errors dues to the thickness of the tubes, and are of order at most $\delta$.   The first assertion is easy:
 \begin{align*}
 |\omega_j-\omega_i|&=|2s(I+sC)(y_j-y_i)-\eta|\\
 &\leq C|y_j-y_i|+|\eta|\\
 &\lesssim \frac{1}{2^{l}}
 \end{align*}
 For the second, begin by eliminating the $\omega$s:
 \begin{equation}
 \label{noomegas}
 (t_j-s)\big(I+(t_j+s)C\big)y_j+\eta = (t_i-s)\big(I+(t_i+s)C\big)y_i+\eps
 \end{equation}
 This can be rearranged to give $y_i$ in terms of $y_j$.  Substituting back into the 2nd of our original equations gives
 \begin{align*}
 \omega_i &=\frac{t_i(t_j-s)}{t_i-s}(I+t_iC)\big(I+(t_i+s)C\big)^{-1}\big(I+(t_j+s)C\big)y_j\\
 &\qquad +\frac{t_i}{t_i-s}(I+t_iC)\big(I+(t_i+s)C\big)^{-1}(\eta-\eps)+\eps\\
 &=\frac{t_i(t_j-s)}{t_i-s}(I+t_jC)y_j+\frac{t_i}{t_i-s}\left(I-sC\right)(\eta-\eps)+\eps
 \end{align*}
 where we have used the fact that $C^2=0$.  Looking back at
 \eqref{nonfatomega} we discover that the first term belongs to the
 intersection of the surface determined by the first two curves with
 the horizontal plane.  So the distance we are interested in is at
 most the absolute value of the other two terms, so at most $\tfrac{C}{|t_i-s|}\delta+\delta$.  Finally we just ensure
 that $|t_i-s|$ is comparable to $|t_j-s|$. From \eqref{noomegas}
 using the fact that on $\supp(a)$ the eigenvalues of $I+x_nC$ are bounded above and
 below, we get
 \begin{align*}
 |t_i-s|2^{-k} &\gtrsim \delta2^{l+m}2^{-k}-\delta\\
 |t_i-s|&\gtrsim\delta2^{l+m}-\delta 2^k\\
 &\gtrsim \delta2^{l+m}
 \end{align*}
 provided that $k-l+m$ is not too large. Since $l\geq k-2$ this could
 happen only with $l$ close to $k$ and $m$ small, in which case the
 claim is trivial anyway.  So the distance of $\omega_i$ from the
 curve of intersection is at most $2^{-(l+m)}$ and we have proved the
 claim.\medskip

 We can now complete the proof of the Hairbrush Lemma, and hence the whole theorem.
 The claim tells us that 
 $$\#H_{j,k,l,m}\lesssim 2^{-l}\left(2^{-(l+m)}\right)^{n-2}\delta^{n-1}$$ which we use as follows:
 \begin{align*}
 M_{H_{j,k,l,m}}&(i) :=\frac{1}{|T_i|}\int_{T_i}\bigg(\sum_{p\in H_{j,k,l,m}}\1_{T_p}(x)\bigg)^\frac{1}{n-1}\d x\\
 &= \frac{1}{|T_i|}\int_{\{x\in T_i:\dist(x,T_i\cap T)\sim\delta2^{l+m}\}}\!\!\bigg(\sum_{p\in H_{j,k,l,m}}\1_{T_p}(x)\bigg)^\frac{1}{n-1}\!\!\d x &&&\mbox{\hspace{-5cm}by defn.~of $H_{j,k,l,m}$}\\
 &\leq \frac{1}{|T_i|}\bigg|\left\{x\in T_i:\dist(x,T_i\cap T)\sim\delta2^{l+m}\right\}\bigg|^{1-\frac{1}{n-1}}
\bigg(\int_{T_i}\sum_{p\in H_{j,k,l,m}}\1_{T_p}(x)\bigg)^\frac{1}{n-1}\d x\\* &&&&\mbox{\hspace{-5cm}by H\"older}\\
 & \leq \frac{1}{\delta^{n-1}}(\delta^n 2^{l+m})^{1-\frac{1}{n-1}}\bigg(\sum_{p\in H_{j,k,l,m}}|T_i\cap T_p|\bigg)^\frac{1}{n-1}\\
 &\lesssim \frac{1}{\delta^{n-1}}(\delta^n 2^{l+m})^{1-\frac{1}{n-1}}\left(\#H_{j,k,l,m}\delta^n2^l\right)^\frac{1}{n-1}&&&\mbox{\hspace{-5cm}by Lemma~\ref{lem:intersec}}\\
 &\lesssim \frac{1}{\delta^{n-1}}(\delta^n 2^{l+m})^{1-\frac{1}{n-1}}\left(2^{-l}\left(2^{-(l+m)}\right)^{n-2}\delta^{n-1}\delta^n2^l\right)^\frac{1}{n-1}&&&\mbox{\hspace{-5cm}by the claim}\\
 &=1.
 \end{align*}
 Summing over all the index sets gives the result.
\end{proof}
\section{Arithmetic Methods}
\label{sec:arith}
\subsection{Introduction}
Sections~\ref{sec:trivial} and \ref{sec:geometric} showed how
geometric methods could give lower bounds for the set dimension of the
form $\frac n2 + \mathrm{const}$.  The best known results in the straight line case in low dimensions ($n=3$ or $4$) are still of this form \cite{katzlabatao:r3,labatao:medium}, but in higher dimensions far better
results are obtained by an arithmetic approach, since these improve
the coefficient of $n$ to something greater than $1/2$.

The arithmetic arises in the form of {\em sumset inequalities}.  For
these we require some notation.
\begin{notation}  Let $\A,\B\subseteq\Z^{n-1}$ be finite sets and let $\G\subseteq\A\times\B$.  For any $(n-1)\times(n-1)$ real matrix $X$ define the $X$-sumset of $\A$ and $\B$ by
$$\A+X\B:=\{a+Xb:(a,b)\in\G\}.$$ In the case $X=-I$ write $\A-\B$
and call it the difference set.
\end{notation}
The structure of sumsets, and inequalities regarding the relative
sizes of sum and difference sets, have been extensively studied by
combinatorialists when the matrix $X$ is an integer multiple of the
identity, but they have generally considered only $\G=\A\times\B$.  See \cite{nathanson:sumset,ruzsa:sumsets}.
 The link with the Kakeya problem was noticed in 1999 by Bourgain
 \cite{bourgain:dimkakeya}, and since then many inequalities with
 $\G\subseteq\A\times\B$ have been proved.  However, the case where
 $X$ is not a multiple of $I$ arises only with curves, and seems to be
 a new problem.

The most general matrix sumset problem with $N+2$ ``slices'' is as follows:
\begin{problem}
\label{prob:matsum}
Let $X_1,X_2,\dots,X_N$ be $(n-1)\times (n-1)$ real matrices. 
To avoid trivialities assume that they are non-zero, distinct, and not equal to $-I$.  Does there exist an $\eps>0$ depending only on the $X_j$s such that for all $\A,\B\in\Z^{(n-1)}$, $\G\subseteq\A\times\B$ we have
$$\#(\A-\B)\leq \max\Big\{\#\A,\#\B,\max_j\#(\A+{X_j}\B)\Big\}^{2-\eps}\,\mbox{?}$$
If so, what is the largest possible $\eps$?
\end{problem}
The idea is to let $\A$ and $\B$ correspond to two horizontal slices through our $\delta$-discretised Kakeya set, and $\G\subseteq\A\times\B$ the set of all pairs which are joined by a line.  Then the difference set $\A-\B$ corresponds to the set of directions, so must be large.  However, if our set had small dimension then $\A$ and $\B$ must be small, and moreover the set of midpoints of the lines, which (for straight lines) has cardinality $\#(\A+\B)$, must also be small.  Inequalities regarding the relative sizes of sumsets and difference sets thus lead to lower bounds for the dimension of Kakeya sets.

In the curved case, we must discover how to determine from two endpoints the location of any other point on the curve, and its direction.
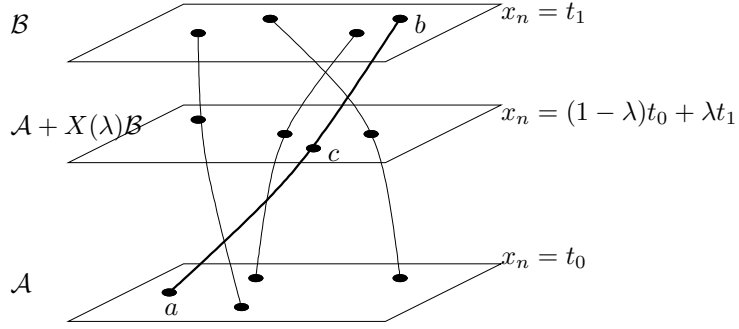
\begin{figure}[!hbtp]
\centering
\setlength{\unitlength}{0.0005in}
{\renewcommand{\dashlinestretch}{30}
\begin{picture}(6469,3342)(0,-10)
\put(2700,3162){\blacken\ellipse{150}{74}}
\put(2700,3162){\ellipse{150}{74}}
\put(4050,462){\blacken\ellipse{150}{74}}
\put(4050,462){\ellipse{150}{74}}
\put(3750,1962){\blacken\ellipse{150}{74}}
\put(3750,1962){\ellipse{150}{74}}
\path(4050,462)(4050,464)(4049,468)
	(4048,476)(4047,489)(4045,507)
	(4042,531)(4038,560)(4034,596)
	(4029,638)(4023,685)(4016,737)
	(4009,792)(4001,851)(3992,912)
	(3983,975)(3974,1038)(3964,1101)
	(3954,1164)(3944,1225)(3934,1285)
	(3924,1343)(3914,1398)(3903,1452)
	(3892,1503)(3882,1552)(3871,1599)
	(3859,1644)(3848,1687)(3836,1728)
	(3824,1767)(3812,1805)(3799,1842)
	(3785,1877)(3771,1912)(3756,1946)
	(3741,1979)(3725,2012)(3708,2044)
	(3691,2076)(3673,2108)(3654,2140)
	(3634,2172)(3612,2204)(3590,2236)
	(3566,2269)(3540,2302)(3513,2336)
	(3485,2371)(3455,2407)(3422,2443)
	(3388,2481)(3353,2521)(3315,2561)
	(3276,2602)(3235,2644)(3193,2687)
	(3150,2731)(3106,2774)(3062,2818)
	(3018,2860)(2975,2902)(2934,2942)
	(2895,2979)(2858,3014)(2825,3045)
	(2795,3073)(2770,3097)(2749,3117)
	(2732,3132)(2719,3144)(2710,3153)
	(2704,3158)(2701,3161)(2700,3162)
\put(1950,3012){\blacken\ellipse{150}{74}}
\put(1950,3012){\ellipse{150}{74}}
\put(1650,312){\blacken\ellipse{150}{74}}
\put(1650,312){\ellipse{150}{74}}
\put(3600,3012){\blacken\ellipse{150}{74}}
\put(3600,3012){\ellipse{150}{74}}
\put(4050,3162){\blacken\ellipse{150}{74}}
\put(4050,3162){\ellipse{150}{74}}
\put(2550,462){\blacken\ellipse{150}{74}}
\put(2550,462){\ellipse{150}{74}}
\put(2400,162){\blacken\ellipse{150}{74}}
\put(2400,162){\ellipse{150}{74}}
\put(1950,2112){\blacken\ellipse{150}{74}}
\put(1950,2112){\ellipse{150}{74}}
\put(3150,1812){\blacken\ellipse{150}{74}}
\put(3150,1812){\ellipse{150}{74}}
\put(2850,1962){\blacken\ellipse{150}{74}}
\put(2850,1962){\ellipse{150}{74}}
\path(1800,3312)(5100,3312)(3900,2712)
	(600,2712)(1800,3312)(1800,3312)
\path(1800,612)(5100,612)(3900,12)
	(600,12)(1800,612)(1800,612)
\path(1800,2262)(5100,2262)(3900,1662)
	(600,1662)(1800,2262)(1800,2262)
\path(2400,162)(2400,164)(2399,167)
	(2397,174)(2394,185)(2390,201)
	(2384,222)(2377,248)(2369,280)
	(2359,317)(2348,359)(2336,406)
	(2322,457)(2308,511)(2294,568)
	(2278,627)(2263,687)(2247,748)
	(2232,808)(2217,868)(2202,927)
	(2188,984)(2174,1040)(2161,1094)
	(2148,1146)(2136,1196)(2125,1245)
	(2114,1291)(2104,1335)(2094,1377)
	(2085,1418)(2077,1458)(2069,1496)
	(2062,1533)(2055,1568)(2048,1603)
	(2042,1637)(2036,1671)(2030,1704)
	(2025,1737)(2019,1775)(2014,1813)
	(2008,1851)(2004,1889)(1999,1927)
	(1995,1966)(1991,2006)(1987,2046)
	(1984,2088)(1980,2131)(1977,2175)
	(1975,2221)(1972,2269)(1969,2319)
	(1967,2370)(1965,2422)(1963,2476)
	(1961,2531)(1959,2585)(1958,2640)
	(1956,2693)(1955,2744)(1954,2793)
	(1953,2838)(1952,2878)(1952,2913)
	(1951,2943)(1951,2966)(1950,2985)
	(1950,2997)(1950,3006)(1950,3010)(1950,3012)
\path(2550,462)(2550,464)(2551,468)
	(2552,475)(2553,486)(2555,503)
	(2558,525)(2561,552)(2565,584)
	(2570,622)(2575,665)(2582,712)
	(2588,763)(2596,817)(2603,873)
	(2611,930)(2620,988)(2628,1046)
	(2637,1103)(2646,1159)(2655,1214)
	(2664,1267)(2672,1318)(2681,1367)
	(2690,1414)(2700,1459)(2709,1502)
	(2718,1544)(2727,1583)(2737,1621)
	(2747,1658)(2757,1693)(2767,1727)
	(2778,1761)(2789,1793)(2801,1825)
	(2813,1856)(2825,1887)(2839,1919)
	(2853,1951)(2868,1983)(2883,2015)
	(2900,2047)(2917,2080)(2936,2113)
	(2956,2146)(2977,2181)(2999,2216)
	(3023,2252)(3048,2290)(3074,2329)
	(3103,2369)(3132,2410)(3163,2453)
	(3196,2497)(3229,2541)(3264,2587)
	(3298,2632)(3333,2678)(3368,2722)
	(3402,2765)(3434,2806)(3465,2844)
	(3493,2879)(3518,2910)(3539,2937)
	(3557,2960)(3572,2978)(3583,2991)
	(3591,3001)(3596,3007)(3599,3011)(3600,3012)
\thicklines
\path(1650,312)(1651,313)(1654,315)
	(1659,320)(1667,327)(1678,337)
	(1693,351)(1712,368)(1735,389)
	(1762,414)(1793,443)(1828,474)
	(1866,509)(1907,547)(1950,587)
	(1995,628)(2041,671)(2088,715)
	(2135,759)(2183,804)(2230,848)
	(2276,891)(2321,934)(2365,976)
	(2408,1017)(2449,1057)(2488,1096)
	(2527,1134)(2564,1170)(2599,1206)
	(2633,1240)(2666,1274)(2698,1307)
	(2729,1339)(2758,1371)(2787,1402)
	(2816,1433)(2843,1464)(2870,1494)
	(2897,1525)(2924,1556)(2950,1587)
	(2976,1618)(3002,1650)(3028,1682)
	(3054,1714)(3080,1747)(3107,1781)
	(3133,1816)(3160,1852)(3188,1889)
	(3216,1927)(3245,1967)(3275,2008)
	(3305,2051)(3337,2095)(3369,2142)
	(3403,2190)(3437,2240)(3472,2292)
	(3509,2345)(3546,2400)(3584,2456)
	(3622,2513)(3661,2570)(3699,2628)
	(3737,2686)(3774,2742)(3811,2797)
	(3845,2849)(3878,2899)(3908,2945)
	(3936,2988)(3961,3025)(3983,3058)
	(4001,3087)(4016,3110)(4028,3128)
	(4037,3142)(4043,3151)(4047,3157)
	(4049,3161)(4050,3162)
\put(1600,82){$a$}
\put(4200,3012){$b$}
\put(3300,1682){$c$}
\put(0,312){$\A$}
\put(0,1962){$\A+X(\lambda)\B$}
\put(0,3012){$\B$}
\put(5100,3162){$x_n=t_1$}
\put(5100,2112){$x_n=(1-\lambda)t_0+\lambda t_1$}
\put(5100,612){$x_n=t_0$}
\end{picture}
}
\caption{Slices through a curved Kakeya set}\label{fig:curveslice}
\end{figure}
\begin{lem}
\label{lem:derivemat}
Let $E$ be a set of curves of the form \eqref{curves} and let $\A,\B\subset\Z^{n-1}$ be the ($\delta$-discretised) intersections of $\nbd(E)$ with the planes $x_n=t_0$ and $x_n=t_1$ respectively, where $t_0\neq t_1$. Let $$\G:=\{(a,b): a\mbox{ and }b\mbox{ lie on the same tube in }\nbd(E)\}\subseteq\A\times\B.$$  Then 
\begin{enumerate}
\item The set of directions $y$ has the same cardinality as the difference set $\A-\B$.
\item Assume that $t_0,t_1\neq0$.  The set of centres $\omega$ has the same cardinality as $\A-T\B$ where $T$ is the $(n-1)\times(n-1)$ matrix
$$T=\frac{t_0}{t_1}(I+t_0C)(I+t_1C)^{-1}.$$
\item
The intersection of the set with the plane $x_n= (1-\lambda)t_0+\lambda t_1$ has the same cardinality as the sumset $\A+X(\lambda)\B$, where $X(\lambda)$ is the $(n-1)\times(n-1)$ matrix
\begin{multline*}
\hspace*{-0.07\textwidth}X(\lambda)=\frac{\lambda}{1-\lambda}\left[I+\lambda(t_1-t_0)C\big(I+(t_0+t_1)C\big)^{-1}\right]^{-1}\Big[I-
\left.(1-\lambda)(t_1-t_0)C\big(I+(t_0+t_1)C\big)^{-1}\right].
\end{multline*}
\end{enumerate}
\end{lem}
\begin{proof}
Consider a curve through the points $(a, t_0)$ and $(b,t_1)$.
The equation \eqref{curves} of the curves gives 
\begin{align}
a&=\omega -t_0y-t_0^2Cy\label{a}\\
b&=\omega -t_1y-t_1^2Cy.\label{b}
\end{align} 
Subtracting these we find that
$y=\frac1{t_0-t_1}\big(I+(t_0+t_1)C\big)^{-1}(b-a)$, and so the
first assertion follows, since multiplication by an invertible matrix
does not change the cardinality.  This part is analogous to Lemma~\ref{lem:intersec}.

If we write $M_j=t_j(I+t_jC)$ for $j=0,1$ so that $a=\omega-M_0y$ and $b=\omega-M_1y$, then solving gives
$$\omega=(M_0^{-1}-M_1^{-1})^{-1}(M_0^{-1}a-M_1^{-1}b).$$  
We can always multiply through by an invertible matrix to get the $a$ on its own.  Therefore an appropriate ``difference set'' is $\A-T\B$ where $T=M_0M_1^{-1}=\frac{t_0}{t_1}(I+t_0C)(I+t_1C)^{-1}$ as in the second assertion.  Note that for the Nikodym problem, we cannot take slices through $x_n=0$ anyway, because these sets arise only when $x_n=0$ is not in the support of the cutoff function in \eqref{TN}.

For the third, denote the point of intersection of this curve with the intermediate plane by $c$.  It helps to take $(1-\lambda)\eqref{a}+\lambda\eqref{b}$, which gives
$$\omega = (1-\lambda)a+\lambda b +(1-\lambda)\big(t_0y+t_0^2Cy\big)+\lambda\big(t_1y+t_1^2Cy\big).$$
This allows lots of cancellation, so that
\begin{align*}
c &= \omega-\big((1-\lambda)t_0+\lambda t_1\big)y-\big((1-\lambda)t_0+\lambda t_1\big)^2 Cy\\
 &= (1-\lambda)a+\lambda b +\lambda(1-\lambda)(t_0-t_1)^2Cy\\
 &= (1-\lambda)a+\lambda b +\lambda(1-\lambda)(t_0-t_1)^2C\frac1{t_0-t_1}\big(I+(t_0+t_1)C\big)^{-1}(b-a)\\
 &= \left[(1-\lambda)I+\lambda(1-\lambda)(t_1-t_0)C\big(I+(t_0+t_1)C\big)^{-1}\right]a+\\& \hspace{0.3\textwidth}+\left[\lambda I-\lambda(1-\lambda)(t_1-t_0)C\big(I+(t_0+t_1)C\big)^{-1}\right]b.
\end{align*}
Multiplying through by an invertible matrix gives the result.
\end{proof}
It is easy to check that all the matrices occurring above are indeed
 invertible, because of the non-degeneracy criterion \eqref{nondegAB}.
 Recall also that in the straight line case we have $C=0$ and hence
 $X(\lambda)$ is really just a scalar.  The sumset in the second
 assertion does not appear in the literature on the straight line
 problem since it is only appropriate when dealing with Nikodym rather than
 Kakeya sets.  Although the elements of the matrices $C$ and the
 parameters $t_0,t_1,\lambda$ are real, for our application we should
 consider only matrices over $\Q$, since each real number may be
 approximated to within $O(\delta)$ by a rational, which corresponds
 to the same point in the $\delta$-discretisation.

We now show how sumset inequalities imply results about Kakeya sets.  For straight line Kakeya sets this is due to Bourgain \cite{bourgain:dimkakeya}, with the ``plane-varying'' improvement noticed by Katz and Tao \cite{katztao:sumsets}.
\begin{lem}
\label{lem:application}
\begin{enumerate}
\item
Suppose that for some matrix $C$ we can choose
$t_0,t_1\in[-1,1]$ and $\lambda_j\in(0,1),\ j=1,\dots,N$ such that
Question~\ref{prob:matsum} with $X_j=X(\lambda_j)$ has a positive answer.  Then Kakeya sets of curves of the form \eqref{curves} for
this $C$ have Minkowski dimension at least $\frac{n-1}{2-\eps}$.
\item If the same holds but with $X_j=X(\lambda_j)T^{-1}$ where $T$ is as in Lemma~\ref{lem:derivemat} and none of the heights $t_0,t_1,(1-\lambda_j)t_0+\lambda_jt_1$ is $0$, then the corresponding curved Nikodym sets have Minkowski dimension at least $\frac{n-1}{2-\eps}$.
\end{enumerate}
\noindent In both cases, if in fact we have a range of solutions, meaning that Question~\ref{prob:matsum} remains true as $t_0$ is allowed to vary over some small interval and the other heights to vary correspondingly, then we can obtain the better lower bound of $\frac{n-1}{2-\eps}+1$.
\end{lem}
\begin{proof}
Let $E$, $\A$, $\B$ and $\G$ be as in Lemma~\ref{lem:derivemat}.
\begin{proofenum}
\item If $E$ is a curved Kakeya set, then we may assume that $\nbd(E)$ consists of $\delta^{-(n-1)}$ tubes in distinct $\delta$-separated directions.  So by Lemma~\ref{lem:derivemat} we have $\delta^{-(n-1)}\sim\#(\A-\B)$.  By the assumption, this means that
$$\delta^{-(n-1)}\lesssim
\max\Big\{\#\A,\#\B,\max_j\#(\A+{X(\lambda_j)}\B)\Big\}^{2-\eps}$$
and so one of the sets on the right hand side has cardinality at least $\delta^{-\frac{n-1}{2-\eps}}$.  So $\nbd(E)$ includes a $\delta$-ball at each of these points, so that $|\nbd(E)|\gtrsim \delta^n\delta^{-\frac{n-1}{2-\eps}}$, which says that $E$ has Minkowski dimension at least $\frac{n-1}{2-\eps}$.
\item If $E$ is a curved Nikodym set, then we may assume that $\nbd(E)$ consists of $\delta^{-(n-1)}$ tubes whose centres $\omega$ are distinct and $\delta$-separated.  So by Lemma~\ref{lem:derivemat}, we have $\delta^{-(n-1)}\sim\#(\A-T\B):=\#(\A-\B')$ where $\B':=\{Tb:b\in\B\}$.  By the assumption, this means that
$$\delta^{-(n-1)}\lesssim
\max\Big\{\#\A,\#\B',\max_j\#(\A+{X(\lambda_j)}T^{-1}\B')\Big\}^{2-\eps}$$
and so one of the sets on the right hand side has cardinality at least $\delta^{-\frac{n-1}{2-\eps}}$, and since $T^{-1}\B'=\B$, this implies that $\nbd(E)$ includes  $\delta$-balls as before, giving the same bound for the dimension.
\end{proofenum}
If we have a range of solutions, as we vary the heights $t_0,t_1$ of $\A,\ \B$
and hence those of  the other $N$ slices, we always have at least one of the intersections
having large cardinality.  So one of them has this
property for a range of heights of positive measure.  This means that
instead of $\delta$-balls, $\nbd(E)$ actually includes small cylinders of width $\delta$ and height some constant $c$.  Hence
$$|\nbd(E)| \gtrsim c\delta^{n-1-\frac{n-1}{2-\eps}}$$
which shows that $\dim(E)\geq\frac{n-1}{2-\eps}+1$ as required.
\end{proof}
So we must try to answer Question~\ref{prob:matsum}. Clearly it always holds with $\eps=0$; this recovers the ``trivial bound'' for which we proved the stronger maximal function version in Section~\ref{sec:trivial}.  If we could obtain $\eps=1$ then the sets would have dimension $n$; unfortunately this is not known for any matrices, and in the scalar case with $N=1$ has been shown to be false by Ruzsa \cite{ruzsa:sumsets}).
\subsection{The scalar case}
We will begin by reviewing the known results in the case where all of the
$X_j$ are multiples of the identity, and seeing what results for curves can be deduced from them.

With three slices ($N=1$) and $X_1=I$ we have the problem Bourgain
originally used in \cite{bourgain:dimkakeya}.  He proved the estimate
with $\eps=\frac1{13}$, which was quickly improved to $\frac16$ by
Katz and Tao \cite{katztao:sumsets}, i.e.
\begin{equation}\label{sixth}
\#(\A-\B)\leq \max\big\{\#\A,\#\B,\#(\A+\B)\big\}^{2-1/6}.\end{equation}
  For all other rational multiples
of $I$ the existence of positive improvements $\eps$ has been proved
by Christ \cite{christ:trilinear}, although it is tedious to compute
their values.

The first four-slice estimate was again due to Katz and Tao, namely
\begin{equation}
\label{quarter}
\#(\A-\B)\leq \max\big\{\#\A,\#\B,\#(\A+\B),
\#(\A+2\B)\big\}^{2-1/4}\end{equation} as shown in
\cite{katztao:sumsets}.  In \cite[Theorem 3.3]{katztao:newbounds} they
showed that $\eps=1/4$ still holds if, instead of using $1$ and
$2$ as here, the two non-zero scalars simply differ by $1$.
We shall generalise this for matrices shortly.

In the same theorem, they showed that six slices with scalars
$x,y,\bar{x},\bar{y}$ satisfying
\begin{equation}
\label{relation}
(1+\tfrac1{\bar{x}})x=(1+\tfrac1{\bar{y}})y\end{equation}
 also gave $\eps=1/4$.  This relation allows us to obtain results for five slices also, by taking two scalars to be equal (or by taking one to be $\infty$, in which case we interpret $\A+\infty\B=:\B$).

They also proved an iteration result:
\begin{thm}[Katz \& Tao \cite{katztao:newbounds}]\label{thm:iter}
If we can obtain $\eps=\eps_0$ in Question~\ref{prob:matsum} for some
finite set of scalars, then for some larger set of scalars we can
obtain $\eps=\frac{2-\eps_0^2}{8-7\eps_0+\eps_0^2}$.  Hence by
choosing larger and larger sets, the improvement $\eps$ may be made as
close to the fixed point $0.32486\dots$ as we wish.
\end{thm}
This result gives the lower bound of approximately\label{approx}
$0.5969n+0.403$ for the Minkowski dimension of straight-line Kakeya
sets, which is currently the best known for large $n$. 

In order to apply these results in the curved case we must discover whether $X(\lambda)$, or $X(\lambda)T^{-1}$ can be multiples of the identity.
\begin{lem} \begin{enumerate}
\item $X(\lambda)$ cannot be a multiple of the identity except in the straight line case.
\item $X(\lambda)T^{-1}$ is a multiple of the identity for all $\lambda$ if $C^2=0$, but this condition is not necessary to obtain the equality for some $\lambda$.
\end{enumerate}
\end{lem}
\begin{proof}
It will be helpful to write \label{defc}$M=M_{t_0,t_1}:=(t_1-t_0)C\big(I+(t_0+t_1)C\big)^{-1}$ so that $X(\lambda):=\frac{\lambda}{1-\lambda}[I+\lambda M]^{-1}[I-(1-\lambda)M]$.
Suppose that $X(\lambda)=\frac{\lambda}{1-\lambda}kI$ where $k=k(t_0,t_1,\lambda)$ is some scalar function.  Then
\begin{align*}
I-(1-\lambda) M &=k(I+\lambda M)\\
(1-k)I &= (k\lambda+1-\lambda)M
\end{align*}
which implies that $M$ is some (possibly zero) multiple of $I$.  By the definition of $M$ this implies that
$C$ is a multiple of $I$.  As observed before (page \pageref{straighten}) this reduces to the straight line case.

On the other hand, if $C^2=0$, then $M=(t_1-t_0)C$ and hence $X(\lambda)=\frac{\lambda}{1-\lambda}(I-(t_1-t_0)C)$, while $T:=\frac{t_0}{t_1}(I+t_0C)(I+t_1C)^{-1} = \frac{t_0}{t_1}(I-(t_1-t_0)C)$. So $X(\lambda)$ and $T$ are parallel. 
Theorem~\ref{thm:nastynik} gives examples where $X(\lambda)=T$ but $C^2\neq0$.
\end{proof}
\begin{thm}[Nikodym result for $C^2=0$]\label{thm:bestnik}
Nikodym sets of curves of the form \eqref{curves} with $C^2=0$ have Minkowski dimension at least $\frac{n-1}{2-\eps}+1\approx0.5969n+0.403$, where $\eps$ is the smallest root of $\eps^3-6\eps^2+8\eps-2$.
\end{thm}
\begin{proof}In this case, for all $t_0,t_1$ the sumsets are just the scalar ones $\A+\frac{\lambda}{1-\lambda}\B$.  Clearly by choosing suitable heights these can be any scalars we like, so this follows immediately from Katz and Tao's sumset result (Theorem~\ref{thm:iter}).
\end{proof} 
Many other families of curves admit some good bound for the Nikodym sets, however.  To use Katz and Tao's simple three-slice estimate \eqref{sixth}, all we require is that there exist $t_0,t_1\in[-1,1]\setminus\{0\}$ and $\lambda\in(0,1)\setminus\{\frac{t_0}{t_0-t_1}\}$ such that $X(\lambda)=T$.  These equations are difficult to solve, but where the matrix $C$ is invertible, or where one of $A,B$ is positive definite so that we may assume that $C$ is diagonal and hence commuting, we can simplify the problem.
\begin{thm}\label{thm:nastynik} Suppose that $C$ is either diagonal or invertible.
Then $X(\lambda)=T$ if and only if $C$ satisfies the quadratic equation
\begin{equation}\label{nastynikquad}\left(t_0^2t_2^2+t_1^2t_2^2-2t_0^2t_1^2\right)C^2+(t_0+t_1+t_2)(t_0t_2+t_1t_2-2t_0t_1)C+(t_0t_2+t_1t_2-2t_0t_1)I=0\end{equation}
where  $t_0,t_1$ and $t_2:=
(1-\lambda)t_0+\lambda t_1$ are the heights of the planes.  If this is so for some choice of $t_0,t_1,t_2\in (-1,1)$,  then the corresponding curved Nikodym sets have Minkowski dimension at least $\frac{6n-6}{11}$, while if there is a whole range of such heights then the dimension is at least $\frac{6n+5}{11}$. 
\end{thm}\begin{rmk}
For any such heights to exist, $C$ must have at most two eigenvalues $h$ and $k$ with $|\frac1h+\frac1k|<3$.  Some cases in which $h$ and $k$ satisfy this quadratic for a range of heights are as follows:
\begin{enumerate}
\item If he eigenvalues $h$ and $k$ are real, then the necessary condition $|\frac1h+\frac1k|<3$ is also sufficient.
\item If they are complex, so $h=\bar{k}=\alpha+i\beta$, $\alpha,\beta\in\R$, then it is enough for either of the following to hold: 
\begin{enumerate}
\item $\max\left\{2|\alpha|-\alpha^2,\frac12\left(-1-2\alpha^2+\sqrt{1+16\alpha^2}\right)\right\}\leq\beta^2\leq3\alpha^2$.
\item $\beta^2\geq\max\left\{2|\alpha|-\alpha^2,\frac12\left(-1-2\alpha^2+\frac12\sqrt{16\alpha^2+24|\alpha|+1}\right)\right\}$.
\end{enumerate}
In particular, $|\alpha|\geq2.36\dots$ suffices.
\end{enumerate}
\end{rmk}
\begin{proof}[Proof of theorem and remark]
When $C$ is diagonal or invertible, the equation $X(\lambda)=T$ can be easily rearranged to give \eqref{nastynikquad}.  Thus for these cases, we require that $C$ should have at most two eigenvalues.  Moreover, if the eigenvalues are distinct, then we require $C$ to be diagonalisable (over $\Complex$), while in the case of one repeated eigenvalue we need the Jordan normal form of $C$ to contain only $1\times1$ and $2\times2$ blocks. By considering the ratio of the last two coefficients we find that the sum of the reciprocals of the eigenvalues must be  equal\label{sec:sumrecip} to $-(t_0+t_1+t_2)\in(-3,3)$, which imposes further restriction on $C$.  
\begin{proofenum}
\item
If $C$ has real eigenvalues, then they must lie in $(-1/2,1/2)$ by \eqref{nondegAB}, whence the necessary condition $|\frac1h+\frac1k|<3$ implies that they are of opposite sign.  Without loss of generality $|h|\leq|k|$.  We exploit homogeneity by setting $t_0=t$, $t_1=bt$, $t_2=ct$.  Then \eqref{nastynikquad} becomes 
\begin{equation} 
Q(X):=(2b^2-b^2c^2-c^2)X^2+(b+c+1)(2b-bc-c)X+(2b-bc-c)=0
\end{equation}
which should have roots $ht$ and $kt$.  These lie in $(-1/2,1/2)$ and have opposite sign so that the sum of their reciprocals is in $(-3,3)$.  So we choose $b,c$ so that $Q(0)>0$, $Q(1/2)<0$ and $Q(-1/2)<0$.  It is easy to check that this is so if we choose $b\in(0,1)$ and 
$$\frac{-6-4b-2b^2+2\sqrt{7b^4+28b^3+52b^2+48b+9}}{2(b^2+2b+3)}<c<\frac{2b}{1+b}.$$
Consider the region of those $b,c$ satisfying this for which $b+c+1>|\frac1h+\frac1k|$.  This is shown in Figure~\ref{fig:region}, and is not empty provided that $|\frac1h+\frac1k|<3$.
\begin{figure}[!htbp]
\begin{center}
\setlength{\unitlength}{0.0006in}
\begin{picture}(4098,3885)(0,-10)
\put(6890.556,-2894.714){\arc{14677.445}{3.5739}{3.9714}}
\put(6942.493,-2940.740){\blacken\arc{14815.484}{3.9705}{4.2783}}
\put(6942.493,-2940.740){\arc{14815.484}{3.9705}{4.2783}}
\blacken\path(1937,2520)(2027,2430)(3827,3780)
	(2927,3240)(1937,2520)
\put(9002.000,-4995.000){\whiten\arc{20374.617}{3.6744}{4.1795}}
\put(9002.000,-4995.000){\arc{20374.617}{3.6744}{4.1795}}
\path(227,180)(3827,180)
\path(3707.000,150.000)(3827.000,180.000)(3707.000,210.000)
\path(227,180)(227,3780)
\path(257.000,3660.000)(227.000,3780.000)(197.000,3660.000)
\path(677,3780)(3827,630)
\put(227,3820){\makebox(0,0)[lb]{$c$}}
\put(3952,180){\makebox(0,0)[b]{$b$}}
\put(3827,0){\makebox(0,0)[b]{$1$}}
\put(47,3600){\makebox(0,0)[b]{$1$}}
\end{picture}
\caption{The region $\frac{-6-4b-2b^2+2\sqrt{7b^4+28b^3+52b^2+48b+9}}{2(b^2+2b+3)}<c<\frac{2b}{1+b}$ and $b+c+1>|\frac1h+\frac1k|$.}
\label{fig:region}
\end{center}
\end{figure}
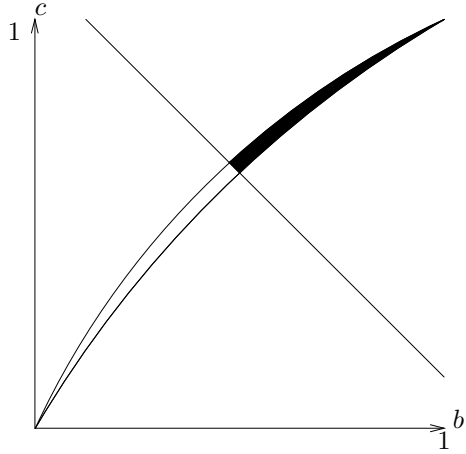

Let $X^+(b,c),X^-(b,c)$ denote the roots of $Q$ obtained by the quadratic formula, taking the positive or negative square root respectively, and 
define the function $f(b,c):=-(b+c+1)X^+(b,c)$, which is continuous on the interior of the region.  

On the upper curve we have $X^+=X^-=0$ and so $f=0$, while on the lower curve $X^-=1/2$.  The ratio of the last two coefficients of the quadratic function $Q$ tells us that $\frac1{X^+}+\frac1{X^-}=-(b+c+1)$, so that on the lower curve we have $f(b,c)=\frac{b+c+1}{b+c+3}$ which tends to $3/5$ as $(b,c)\to(1,1)$.  So provided that $1+h/k\in(0,3/5)$ we can find a curve of points $(b_0,c_0)$ on which $f(b_0,c_0)=1+h/k$. We have $|X^+(b_0,c_0)|=|\frac{1+h/k}{b_0+c_0+1}|<\frac{|1+h/k|}{|\frac1h+\frac1k|}=|h|$, so that the heights $t:=X^+(b_0,c_0)/h$, $t_1:=b_0t$ and $t_2:=c_0t$ are all in $(-1,1)$.  Finally observe that 
$$1+h/k=-X^+(b_0,c_0)(b_0+c_0+1)=X^+\left(\frac1{X^+}+\frac1{X^-}\right)=1+th/X^-$$
so that $X^-(b_0,c_0)=tk$ as required.  Observing that the inequalities $|\frac1h+\frac1k|<3$ and $|h|\leq|k|$ imply $0<1+h/k<3/5$ for $h,k\in(-1/2,1/2)$ gives the result stated.
\item If the eigenvalues of $C$ are complex conjugates $\alpha\pm i\beta$, consider the map $g:\R^3\to\R^2$ given by
$$g(t_0,t_1,t_2)=\left(-(t_0+t_1+t_2),\frac{t_0t_2+t_1t_2-2t_0t_1}{t_0^2t_2^2+t_1^2t_2^2-2t_0^2t_1^2}\right).$$
The two expressions on the right are the sum of the reciprocals of the roots and the product of the roots respectively.  So we have to solve $g(t_0,t_1,t_2)=\big(\frac{2\alpha}{\alpha^2+\beta^2},\alpha^2+\beta^2\big)$.
\begin{enumerate}
\item[(a)] Take $t_1=-t_0$.  Then solving gives
\begin{align*}
t_0=-t_1&=\frac{\sqrt{3\alpha^2-\beta^2}}{\alpha^2+\beta^2} & t_2&=\frac{-2\alpha}{\alpha^2+\beta^2}
\end{align*}
which are in $(-1,1)$ if and only if $\alpha$ and $\beta$ are as claimed.
\item[(b)] Take $t_2=-t_0$.  This is not so easily solved, but $t_1=\frac{-2\alpha}{\alpha^2+\beta^2}$ while $t_0$ satisfies the following cubic:
$$(\alpha^4+2\alpha^2\beta^2+\beta^4)t_1^3+(\beta^2-3\alpha^2)t_1-6\alpha=0.$$
By considering the sign of this cubic at $t_1=-1,\ 0,\ 1$ we can force a sign change in the interval $(-1,1)$ by letting $\alpha$ and $\beta$ satisfy the inequality stated.
\end{enumerate}
In both cases it is easy to check that the Jacobian of $g$ has full rank, so that we can use the implicit function theorem to find not just one solution for $(t_0,t_1,t_2)$ but a whole range.
\end{proofenum}
Once we have found ranges of heights so that the quadratic is satisfied, the result follows from \eqref{sixth} and Lemma~\ref{lem:application}.
\end{proof}
It is likely that there are many other complex pairs $h,k$ which work, however it seems difficult to describe the set of all such pairs concisely.
\subsection{Non-scalar matrices}
\label{sec:genmatsumset}
In the case of curved Kakeya rather than Nikodym sets, the matrix $X(\lambda)$ is never a multiple of $I$, nor is there a matrix $T$ which it might cancel with.  So we have no option but to try to answer
Question~\ref{prob:matsum} with the $X_j$ not  multiples of $I$.
This seems hard.  However, we have some negative results, and
have been able to generalise one of the positive results from the
scalar case.  We begin with a rather trivial observation.
\begin{lem}
\label{lem:block}
If all of the $X_j$ are block diagonal with blocks of the same size,
then a sumset inequality for these $X_j$ implies one for each of the
sets of blocks, with the same $\eps$.
\end{lem}
\begin{proof}
Obvious by letting $\A,\B$ consist of vectors with zeros everywhere
except in the block of interest.
\end{proof}
The converse seems likely to be false --- we would need not only that
``collisions'' often occur in each block of coordinates, but that they
often occur in all coordinates at the same time.

We now reveal the easy but disappointing fact that three slices is
simply not enough in the matrix case.
\begin{prop} 
If $X$ is not a multiple of the identity, then the power of $2$ in
$$\#(\A-\B)\leq \max\big\{\#\A,\#\B,\#(\A+X\B)\big\}^2$$
is best possible.
\end{prop}
\begin{proof}
Choose a vector $v$ that is not an eigenvector of $X$, and let $\B$
consist of $M$ equally spaced points along this direction.  Set
$\A=\{Xb:b\in\B\}$, that is, $M$ equally spaced points along the
direction $Xv$.  Then with $\G=\A\times\B$, clearly $\#(\A+X\B)$ is about $2M$, while since
$v$ and $Xv$ are linearly independent, $\#(\A-\B)$ is about $M^2$.  
\end{proof} 
Similar observations with more slices give another negative result.
\begin{prop}
\label{prop:secular}
Suppose that there exists $v\in\R^n$ such that all the vectors $X_j v$ are
rational multiples of some fixed vector $w$ which is not parallel to
$v$ itself.  (That is, $v$ is a secular vector of each pair of
matrices, but is not an eigenvector.)  Then there can be no positive
answer to Question~\ref{prob:matsum}.
\end{prop}
This theorem is rather weak, but it does at least rule out the case
where the matrices $X_j$ are all multiples of each other but not of
the identity, and combining this with Lemma~\ref{lem:block} gives
further examples.  This makes sense because taking more than three
slices is not really giving much more information.
\medskip\par
\begin{proof}
We have $X_jv=\frac{p_j}{q_j}w$ where $p_j,q_j$ are non-zero coprime integers.  Let $M>\prod_{i=1}^Np_iq_i$ be a large integer, and set
\begin{align*}
\A &= \left\{n\Big(\prod_{i=1}^Np_iq_i\Big)w:n=1,\dots, M\right\}\\
\B &= \left\{n\Big(\prod_{i=1}^Nq_i\Big)v:n=1,\dots, M\right\}.
\end{align*}
Then if $\G=\A\times\B$, we find that
\begin{align*}
\A+X_j\B &= \left\{p_j\Big(\prod_{k\neq j}q_k\Big)\left[q_j\Big(\prod_{i\neq j}p_i\Big) m+n\right]w:m,n=1,\dots, M\right\}\\
\A-\B & = \left\{m\Big(\prod_{j=1}^Np_iq_i\Big)w-n\Big(\prod_{i=1}^Nq_i\Big)v,n=1,\dots, M\right\}
\end{align*}
and we have the combinatorial task of finding their cardinality.
Since $w$ and $v$ are linearly independent it is obvious that
$\#(\A-\B)=M^2$.  To find the cardinality of $\A+X_j\B$, first write
$Q_j:=q_j\prod_{i\neq j}p_i$. We need to know how many distinct
numbers $(mQ_j+n)$ there are, so first let $m<M$ be fixed.  Since
$M>Q_j$, there are $Q_j$ such numbers between $mQ_j+1$ and $(m+1)Q_j$
inclusive, after which the remaining values overlap with those for the
next $m$.  When $m=M$ we just get all the numbers $MQ_j+1$ up to
$MQ_j+M$.  Hence we find that $\#(\A+X_j\B)= Q_j(M-1)+M$.  So $\A$ and
$\B$ and all the sumsets all have cardinality about $M$ while the
difference set is about $M^2$.  So there can be no positive answer to Question~\ref{prob:matsum}.
\end{proof}
This result shows us exactly what goes wrong in Bourgain's ``worst case'' example, and more generally for the parabolas for which we have proved good results in the Nikodym case.
\begin{cor} For a Kakeya set of curves of the form \eqref{curves} with $C^2=0$ we cannot prove any non-trivial bound by sumset methods.
\end{cor}
\begin{proof}
If $C^2=0$, then we can calculate
$$X(\lambda)=\frac{\lambda}{1-\lambda}(I+(t_0-t_1)C).$$ So we
would need a sumset result where the $X_j$ were all multiples of
each other but not of $I$.  But we have already seen in Proposition~\ref{prop:secular} that in such a case no non-trivial estimate can hold.
\end{proof}
It is interesting
that curves that behave well for Nikodym should not do so for
Kakeya.  We shall discuss this in the final section.

So far this picture looks bleak.  However, we can prove one or two results in the
positive direction, analogous to those in the scalar case.  Here we
generalise the four-slice result \eqref{quarter} to the matrix
setting.  For legibility, write $X_1=X,\ X_2=Y$.
\begin{prop}
If $Y-X=I$ and $X$ is invertible, then
$$\#(\A-\B)\leq\max\big\{\#\A,\#\B,\#(\A+X\B),\#(\A+Y\B)\big\}^{7/4}.$$
\end{prop}
\begin{proof}
This is just as in \cite{katztao:sumsets}, so we only give an outline.
Start by discarding elements of $\G$ until $\#(\A-\B)=\#\G$, and
denote the maximum on the right hand side by $M$.  We need to show that
$\#\G\leq M^{7/4}$.

The idea is to count {\em trapezia}: sets of four elements of $\G$
consisting of two ``sides'' whose endpoints have the same value of $a$
while the endpoints of the remaining two sides share values of $a+Yb$
and $b$ respectively.  More precisely, a trapezium is a set
$$\left\{(a_0,b_0),(a_0,b_0'),(a_1,b_1),(a_1,b_1')\right\}\subseteq\G$$
such that $a_0+Yb_0=a_1+Yb_1$ and $b_0'=b_1'$.

First count the number of pairs in $\G$ that share their value of $a$.  This is
\begin{align*}\#\{(a,b),(a,b')\in\G\} & = \sum_{a\in\A}\#\{b:(a,b)\in\G\}^2\\
& \geq \frac{\#\G^2}{\#\A}
\end{align*}
by Cauchy-Schwarz.  A trapezium consists of two such pairs that share their value of $(a+Yb,b')$, so by Cauchy-Schwarz again we find that the
number of trapezia is at least 
$$\frac{(\#\G^2/\#\A)^2}{\#(\A+Y\B)\#\B}\geq\frac{\#\G^4}{M^4}.$$  
But we also have the
following algebraic fact:
$$a_1-b_1'=(I+X^{-1})(a_0+Xb_0)-X^{-1}(a_0+Xb_0')-Yb_1.$$
So, since $\#(\A-\B)=\#\G$, knowing $(a_0+Xb_0),(a_0+Xb_0')$, and $b_1$ is enough to determine $(a_1,b_1')$, and hence the whole trapezium by substituting back.  So the number of trapezia is at most $M^3$, which together with the lower bound of ${\#\G^4}/{M^4}$ gives the result.
\end{proof}
 To apply this to the Kakeya problem with curves, we need conditions on the curves (in terms of $C$) that guarantee the existence of $t_0,t_1\in[-1,1]$ and $\lambda,
\mu\in (0,1)$ such that $X(\lambda)-X(\mu)=I$.
Unfortunately, in many cases this cannot be done.  This is hardly
surprising, since for fixed $C$, we are trying to satisfy
$(n-1)^2$ equations with only four unknowns.  As with the Nikodym sets, we end up trying to make $C$ satisfy a polynomial, whose coefficients now depend on the four heights.  However, in the Kakeya case we are able to deal with the $\lambda$ and $\mu$ first, independently of $t_0$ and $t_1$, which makes the problem easier.
\begin{prop}
If the matrix $M:=(t_1-t_0)C\big(I+(t_0+t_1)C\big)^{-1}$ 
is nilpotent, then $X(\lambda)-X(\mu)$ is never equal to the identity.
\end{prop}
\begin{proof}
 Let $k$ be the highest power of $M$ that is non-zero.  Then $[I+\lambda M]^{-1}=\sum_0^k(-1)^n\lambda^nM^n$, and hence
$$X(\lambda)-X(\mu)=\left(\frac{\lambda}{1-\lambda}-\frac{\mu}{1-\mu}\right)I+\sum_1^k(-1)^n
M^n\left(\frac{\lambda^n}{1-\lambda}-\frac{\mu^n}{1-\mu}\right).$$ If
this equals the identity, then some linear combination of
$I,M,M^2,\dots,M^k$ is zero.  But this cannot happen because the
minimum polynomial of a nilpotent matrix is $x^{k+1}$.
\end{proof}
This cuts down the list of possible candidates for $C$.  Note
in particular that all matrices with $C^2=0$ make $M$ nilpotent
for every choice of $t_0,t_1$.  Obviously invertible matrices $C$ are
not ruled out, and nor are diagonal matrices, and whenever either $A$
or $B$ is positive definite we can fix $C$ to be
diagonal by a change of coordinates.
\begin{prop}
Suppose that $C$ is invertible or diagonal.  Then $X(\lambda)-X(\mu)=I$ if and only if $M:=(t_1-t_0)C\big(I+(t_0+t_1)C\big)^{-1}$ satisfies the following quadratic:
\begin{equation}\label{quadratic}
M^2 + \left(\frac2\mu+\frac1{1-\mu}-\frac1{1-\lambda}\right)M+\left(\frac{1}{\lambda\mu}-\frac{1}{\mu(1-\lambda)}+\frac{1}{\lambda(1-\mu)}\right)I=0
\end{equation}
This quadratic cannot have both its roots real and in $(-1,1)$, but suitable $\lambda,\ \mu\in(0,1)$ can be chosen to give any desired roots $l$ and $m$ such that $l+m<-2(1+\sqrt2)$.
\end{prop}
\begin{proof}
By H\"ormander's criterion \eqref{nondegAB}, $C$ is invertible or diagonal if and only if $M$ is. It is then easy to rearrange the equation $X(\lambda)-X(\mu)=I$ to give \eqref{quadratic}.
So $M$ must either be diagonalisable with at most two distinct eigenvalues, or have one repeated eigenvalue and have Jordan normal form consisting only of $1\times1$ and $2\times2$ blocks. 
Suppose that $M$ has two eigenvalues $l$ and $m$.  These are either
real or form a complex conjugate pair, and so both their sum and their
product are real.
We obtain two simultaneous equations by considering the sum and product of the roots of \eqref{quadratic}.
\begin{align}
\frac{1}{\lambda\mu}-\frac{1}{\mu(1-\lambda)}+\frac{1}{\lambda(1-\mu)}&=lm\label{productofroots}\\
 \frac{2}{\mu}+\frac1{1-\mu}-\frac1{1-\lambda}&=-(l+m)\label{sumofroots}
\end{align}
Now \eqref{sumofroots} is linear in $\lambda$ so we solve it to obtain $$\lambda = 1-\frac{\mu(1-\mu)}{2-\mu+\mu(1-\mu)(l+m)}.$$
Of course this needs to lie in $(0,1)$.  
Tedious calculation shows that in the case $l+m>-2(1+\sqrt{2})$ it does so for all $\mu\in(0,1)$. For $l+m< -2(1+\sqrt{2})$ it does so provided we take
\begin{align}
\label{muint}
\mu&\in\left(0,\tfrac{2-(l+m)-\sqrt{(l+m+2)^2-8}}{2(1-l-m)}\right)&\mbox{or}&
&\mu&\in\left(\tfrac{2-(l+m)+\sqrt{(l+m+2)^2-8}}{2(1-l-m)},1\right).
\end{align}
Next we substitute this expression for $\lambda$ back into
\eqref{productofroots}.  After rearranging we obtain an equation which
is quartic in $\mu$ and quadratic in $l$ and $m$.  With the help of
MAPLE we express it as
\begin{align}
\begin{split}
0 &= [l m (l + m - 1)]\mu^4+\big[2lm+(l+m)\big((l+m)^2-lm-1\big)\big]\mu^3\\
&\qquad\qquad+\big[(l+m)(4-(l+m))-1-2lm\big]\mu^2+4(1-l-m)\mu-4\end{split}\\\begin{split}
 &= -\big[ \mu^2 (1-\mu) (\mu m + 1)\big]l^2
+\big[ \mu (\mu^2  m - \mu m + \mu - 2) (\mu m - \mu + 2)\big]l\\
&\qquad\qquad+(\mu^2 m - \mu m + \mu - 2) (\mu m - \mu + 2)\end{split}\label{l}\\
&=:q(\mu,l,m).\nonumber
\end{align}
Note that this is a real-valued function of $\mu$.  For the second part
we use a na\"\i ve approach via the intermediate value theorem.
Setting $\mu=0$ gives $-4$, while $\mu=1$ gives $-(l+1)(m+1)$.  This is positive only when $l$ and $m$ are real with exactly one being less than $-1$.  Similarly we obtain a positive result by setting $\mu=\frac{-2}{m-1}$ or $\frac{-2}{l-1}$, but these too are only permitted if $l,m$ are real and one is less than $-1$.  However, if we instead substitute in either endpoint from equation \eqref{muint} (which is allowable if
and only if $l+m< -2(1+\sqrt{2})$) we obtain
$$\frac{(l+m)^2\pm(l+m-2)\sqrt{(l+m+2)^2-8}}{2(l+m-1)^2}$$
which by more tedious rearranging is seen to be positive for all $l,m$.

For the first assertion we show that the maximum of the function $q$
  over the region $(\mu,l,m)\in[0,1]\times[-1,1]\times[-1,1]$ is zero, and moreover that this is attained only for $\mu=1$.
For interior maxima we use the version \eqref{l} of the equation as a quadratic
  in $l$. Its stationary point (a maximum) occurs at
$$l=\frac{(\mu^2 m-\mu m+\mu-2)(\mu m -\mu+2)}{2\mu(1-\mu)(\mu m+1)}.$$
Now if $|m|<1$ then the denominator is positive, so that the whole fraction will be less than $-1$ if
$$(\mu^2 m-\mu m+\mu-2)(\mu m -\mu+2)<-2\mu(1-\mu)(\mu m+1)$$ which rearranges to
$$\mu^2(1-\mu)m^2\big(3(1-\mu)^2+1\big)(\mu m+1)>0.$$
So there is no zero of $\frac{\partial
q}{\partial l}$ in the region, except perhaps when $\mu=1$ and $m=-1$.
We find that $q(1,l,-1)\equiv 0$.  Now to check the other boundaries:
\begin{align*}
q(\mu,1,m) &=-m^2\mu^2(1-\mu^2)-2(\mu m+1)(2-\mu^2)\\
q(\mu,-1,m) &=-m^2\mu^2(1-\mu)^2-2(\mu m+1)(1-\mu)\big((1-\mu)^2+1\big)
\end{align*}
Both of these are clearly non-positive, and give zero only at $(1,l,-1)$ as we have already seen, and at $(1,1,m)$.
\end{proof}
Unfortunately, $M$ cannot have real eigenvalues outwith $(-1,1)$.  This follows from the non-degeneracy criterion \eqref{nondegAB} and that fact that $l$ is an eigenvalue of $M$ if and only if
$$\det\left[I+\bigl((1-\tfrac1{l}  )t_0+(1+\tfrac1{l}  )t_1\bigr)C\right]=0.$$
However, if the eigenvalues are complex conjugate we obtain a Kakeya result as follows.
\begin{thm}\label{thm:kak4}Suppose that $C$ is invertible and diagonalisable over $\Complex$ and has only two eigenvalues $\alpha\pm\beta i$.  Then if either $\alpha$ is sufficiently large ($|\alpha|>\frac{1+\sqrt{2}}{2}$ will do) or $\beta$ is large compared to $\alpha$, there is a lower bound of $\frac{4n+3}{7}$ for the curved Kakeya sets associated to $C$.
\end{thm}
\begin{proof}
We have seen that this holds if we can make $l+m<-2(1+\sqrt{2})$.
But $l+m$ is simply twice the real part of the eigenvalues of $M:=(t_1-t_0)C\big(I+(t_0+t_1)C\big)^{-1}$, so if we let the eigenvalues of $C$ be $\alpha\pm i\beta$ as before, 
we require
$$(t_1-t_0)\frac{\alpha +(t_0+t_1)(\alpha^2+\beta^2)}{1+2(t_0+t_1)\alpha+(t_0+t_1)^2(\alpha^2+\beta^2)}<-(1+\sqrt{2}).$$
It helps to write $t_1=1-2\eps, t_0=-1+\eps$, where $\eps<2/3$ may be taken as small as we wish.  The inequality becomes
$$\frac{\alpha(1-\eps\alpha)-\eps\beta^2}{(1-\eps\alpha)^2+\eps^2\beta^2}<-\frac{1+\sqrt{2}}{2-3\eps}$$
Clearly this is satisfied for small $\eps$ and large $\beta$: Choosing
$\eps\lesssim1/|\alpha|$ shows that $\beta\gtrsim|\alpha|+1$ will
work.  Alternatively if $\alpha<-\frac{1+\sqrt{2}}{2}$ then we simply
need to take $\eps$ very small.  If $\alpha>\frac{1+\sqrt{2}}{2}$ we simply swap $t_0$ and $t_1$.  In all of these cases, we have in fact found a whole family of solutions for varying $\eps$ so there is no problem with using the argument about varying the heights of the planes which gave the extra $+1$ for the dimension bound in Lemma~\ref{lem:application}.
\end{proof}
So we get a non-trivial result in some cases, although it is not easy
to give the criteria any geometric interpretation.

However, about the case where $C$ is not invertible, or where
$C$ has more than two eigenvalues or two real ones, we cannot
say anything other than that the above proof will not work.

We have not yet considered using four slices in the Nikodym case.  This is more complicated, because we require $X(\lambda)-X(\mu)=T$ instead of $I$, which means that we cannot write this in terms of $M$ and so we must look at all four variables $t_0,t_1,\lambda,\mu$ together, rather than in two stages as we did above.  By the methods already used, we can show that if $C$ is diagonal or invertible then it must satisfy a cubic equation, and that the reciprocals of the roots (the eigenvalues of $C$) must have the same sum as minus the heights of the slices, as we found in Section~\ref{sec:sumrecip}.  As one would expect, it is difficult to say anything more than that explicitly.

So what hope is there for the use of arithmetic methods?  If we still
want to use only four slices for cases not covered above, then we shall have to prove a new sumset
result, that is, find a more flexible condition than
$X(\lambda)-X(\mu)=I$ which guarantees that the difference set is not
too much larger than the two original sets and their $X(\lambda)$ and
$X(\mu)$ matrix sumsets.  Or we could instead look at using more
slices---the techniques in \cite{katztao:newbounds} suggest that
relations like
\begin{align*}
0 &=X(\lambda)-X(\mu)+X(\nu)^{-1}X(\lambda)\\
\mbox{or }0&=X(\lambda)-X(\mu)+X(\nu)^{-1}X(\lambda)-X(\kappa)^{-1}X(\mu),
\end{align*}
which are analogues of \eqref{relation}, would suffice.  But of course these lead to higher degree polynomials
in more variables which make it harder to compute sufficient
conditions for suitable solutions to exist.

\section{Final remarks}
\label{sec:discuss}
We have now seen two non-trivial positive results for curves of the form \eqref{curves} where $C^2=0$, namely the $\frac{n+2}{2}$ bound for the Nikodym maximal function (Theorem~\ref{thm:katz}) and the bound of approximately $0.5969n+0.403$ for the Minkowski dimension of the Nikodym sets (Theorem~\ref{thm:bestnik}).  This condition on the matrix $C$
has a surprising link with Bourgain's ``worst case'' example for the curved Kakeya problem.
That example had $C=(\begin{smallmatrix}0&0\\1&0\end{smallmatrix})$ which clearly satisfies the condition.  This is no accident; the following converse is also true:
\begin{prop}
If $C^2=0$, then the corresponding Kakeya sets can have dimension as low as $n-\operatorname{rank}(C)$.  In particular, in odd dimensions, the ``trivial'' lower bound of $\frac{n+1}{2}$ can be attained, while in even dimensions there can be sets with dimension at least as low as $\frac{n+2}{2}$.
\end{prop}
\begin{proof}
The curves are as in \eqref{curves} where we may assume that $C$ is in rational canonical form.  Then $C^2=0$ if and only if $C$ consists only of $1\times1$ blocks $(\,0\,)$ and $2\times2$ blocks $(\begin{smallmatrix}0&1\\0&0\end{smallmatrix})$.  The rank of $B$ is the number of $2\times 2$ blocks, which can be at most $\frac{n-1}{2}$ if $n$ is odd, and $\frac{n-2}{2}$ if $n$ is even. 
We can now imitate the proof on page \pageref{pushtosurf} in each block, and the result follows.  
\end{proof}
So it seems that the same curves that allow no good bound in the Kakeya case are particularly amenable to proving good bounds in the Nikodym case. This is rather curious, and may reveal a kind of duality between the Kakeya and Nikodym problems. Up until now Restriction/Kakeya and Bochner-Riesz/Nikodym have been thought of as essentially they same \cite{carbery:paraboloid,tao:brimplies}, but the above suggests they might be better described as dual in some way, or even opposite.  This idea is not so strange when one remembers that curvature of the surface in question is {\em good} when considering Restriction (since it causes decay of the Fourier transform) but {\em bad} for Bochner-Riesz (Bochner-Riesz for squares is trivial).

This also shows the importance of Carbery's transformation $$(x',x_n)\mapsto(x'/x_n,1/x_n)$$ which relates the two classes of problems.  This does not preserve parabolas; rather it maps them to hyperbolas.  Another way of phrasing the above is that if for a given matrix, the parabolas can be tightly packed, then the hyperbolas cannot. This would leave straight lines as an overlapping middle case, the only family that this transformation leaves unchanged.  These ideas will be explored further in \cite{carberyw:3notes}.

\bigskip
\noindent Laura Wisewell,
\\*School of Mathematics,
\\*University of Edinburgh, King's Buildings,
\\*Edinburgh,
\\*EH9 3JZ,
\\*Scotland,
\\*United Kingdom\smallskip\par
\noindent\url{l.wisewell@ed.ac.uk}\\*
\url{http://www.maths.ed.ac.uk/~laura}


\begin{thebibliography}{10}

\bibitem{bourgain:besicovitchtype}
J.~Bourgain.
\newblock Besicovitch type maximal operators and applications to {F}ourier
  analysis.
\newblock {\em Geom.\ Funct.\ Anal.\/}, {\bf 1}: 147--187, 1991.

\bibitem{bourgain:lp}
J.~Bourgain.
\newblock {$L^p$} estimates for oscillatory integrals in several variables.
\newblock {\em Geom.\ Funct.\ Anal.\/}, {\bf 1}: 321--374, 1991.

\bibitem{bourgain:dimkakeya}
J.~Bourgain.
\newblock On the dimension of {K}akeya sets and related maximal inequalities.
\newblock {\em Geom.\ Funct.\ Anal.\/}, {\bf 9}: 256--282, 1999.

\bibitem{carbery:covering}
A.~Carbery.
\newblock Covering lemmas revisited.
\newblock {\em Proc. Edinburgh Math. Soc. (2)\/}, {\bf 31}: 145--150, 1988.

\bibitem{carbery:paraboloid}
A.~Carbery.
\newblock Restriction implies {B}ochner-{R}iesz for paraboloids.
\newblock {\em Math.\ Proc.\ Cambridge Philos.\ Soc.\/}, {\bf 111}: 525--529,
  1992.

\bibitem{carbery:igari}
A.~Carbery, E.~Hern{\'a}ndez, and F.~Soria.
\newblock Estimates for the {K}akeya maximal operator on radial functions in
  ${{\mathbb {R}}}\sp n$.
\newblock In S.~Igari, ed., {\em Harmonic analysis (Sendai, 1990)\/}, pp.
  41--50. Springer, Tokyo, 1991.

\bibitem{carberyw:3notes}
A.~Carbery and L.~Wisewell.
\newblock Three notes on {R}estriction and {B}ochner-{R}iesz.
\newblock In preparation.

\bibitem{christ:trilinear}
M.~Christ.
\newblock On certain elementary trilinear operators.
\newblock {\em Math.\ Res.\ Lett.\/}, {\bf 8}: 43--56, 2001. \url{http://math.berkeley.edu/~mchrist/preprints.html}

\bibitem{cdr:max}
M.~Christ, J.~Duoandikoetxea, and J.~L. Rubio~de Francia.
\newblock Maximal operators related to the {R}adon transform and the
  {C}alder\'on-{Z}ygmund method of rotations.
\newblock {\em Duke Math.~J.\/}, {\bf 53}(1): 189--209, 1986.

\bibitem{falconer:frac}
K.~J. Falconer.
\newblock {\em The geometry of fractal sets\/}.
\newblock No.~85 in Cambridge Tracts in Mathematics. Cambridge University
  Press, 1985.

\bibitem{feffC:ball}
C.~Fefferman.
\newblock The multiplier problem for the ball.
\newblock {\em Ann.\ of Math.\/}, {\bf 94}: 330--336, 1971.

\bibitem{hormander:oscill}
L.~H\"ormander.
\newblock Oscillatory integrals and multipliers on {$FL^p$}.
\newblock {\em Ark. Mat.\/}, {\bf 11}: 1--11, 1973.

\bibitem{katz:socialdensity}
N.~H. Katz.
\newblock Social density and the {K}akeya maximal operator.
\newblock {\em To appear\/}.

\bibitem{katzlabatao:r3}
N.~H. Katz, I.~{\L}aba, and T.~Tao.
\newblock An improved bound on the {M}inkowski dimension of {B}esicovitch sets
  in ${{\mathbb R}}^3$.
\newblock {\em Ann.\ of Math.\/}, {\bf 152}: 383--446, 2000. \url{arXiv:math.CA/9903166}

\bibitem{katztao:sumsets}
N.~H. Katz and T.~Tao.
\newblock Bounds on arithmetic projections, and applications to the {K}akeya
  conjecture.
\newblock {\em Math.\ Res.\ Lett.\/}, {\bf 6}: 625--630, 1999. \url{arXiv:math.CO/9906097}

\bibitem{katztao:newbounds}
N.~H. Katz and T.~Tao.
\newblock New bounds for {K}akeya problems.
\newblock {\em J.~Anal.\ Math.\/}, {\bf 87}: 231--263, 2002. \url{arXiv:math.CA/0102135}

\bibitem{labatao:medium}
I.~{\L}aba and T.~Tao.
\newblock An improved bound for the {M}inkowski dimension of {B}esicovitch sets
  in medium dimension.
\newblock {\em Geom.\ Funct.\ Anal.\/}, {\bf 11}: 773--806, 2001. \url{arXiv:math.CA/0004015}

\bibitem{minicozzi:negative}
W.~P. Minicozzi~II and C.~D. Sogge.
\newblock Negative results for {N}ikodym maximal functions and related
  oscillatory integrals in curved space.
\newblock {\em Math.\ Res.\ Lett.\/}, {\bf 4}: 221--237, 1997. \url{arXiv:math.CA/9912202}

\bibitem{nathanson:sumset}
M.~B. Nathanson.
\newblock {\em Additive Number Theory: Inverse Problems and the Geometry of
  Sumsets\/}.
\newblock No. 165 in Graduate Texts in Mathematics. Springer Verlag, New York,
  1996.

\bibitem{ruzsa:sumsets}
I.~Z. Ruzsa.
\newblock Sums of finite sets.
\newblock In D.~V. Chudnovsky, G.~V. Chudnovsky, and M.~B. Nathanson, eds.,
  {\em Number theory New York Seminar 1991--1995\/}, pp. 281--293. Springer,
  New York, 1996.

\bibitem{schlag:geometricineq}
W.~Schlag.
\newblock A geometric inequality with applications to the {K}akeya problem in
  three dimensions.
\newblock {\em Geom.\ Funct.\ Anal.\/}, {\bf 8}(3): 606--625, 1998.

\bibitem{stein:oscill}
E.~M. Stein.
\newblock Oscillatory integrals in {F}ourier analysis.
\newblock In E.~M. Stein, ed., {\em Beijing lectures in harmonic analysis
  (Beijing, 1984)\/}, no. 112 in Annals of Mathematics Studies, pp. 307--355.
  Princeton University Press, New Jersey, 1986.

\bibitem{tao:brimplies}
T.~Tao.
\newblock The {B}ochner-{R}iesz conjecture implies the restriction conjecture.
\newblock {\em Duke Math.~J.\/}, {\bf 96}: 363--376, 1999. \url{http://www.math.ucla.edu/~tao/preprints/kakeya.html}

\bibitem{wisewell:meas0}
L.~Wisewell.
\newblock Families of surfaces lying in a null set.
\newblock {\em To appear in Mathematika\/}, 2003. \url{arXiv:math.CA/0401249}

\bibitem{wisewell:thesis}
L.~Wisewell.
\newblock {\em Oscillatory Integrals and Curved {K}akeya Sets\/}.
\newblock Ph.D. thesis, University of Edinburgh, 2003.
\newblock \url{http://www.maths.ed.ac.uk/~laura/thesis.html}.

\bibitem{wolff:improvedbound}
T.~Wolff.
\newblock An improved bound for {K}akeya type maximal functions.
\newblock {\em Rev.\ Mat.\ Iberoamericana\/}, {\bf 11}(3): 651--674, 1995.

\bibitem{wolff:recentkakeya}
T.~Wolff.
\newblock Recent work connected with the {K}akeya problem.
\newblock In H.~Rossi, ed., {\em Prospects in Mathematics (Princeton, New
  Jersey, 1996)\/}, pp. 129--162. American Mathematical Society, Providence,
  RI, 1999.

\end{thebibliography}
\end{document}